\documentclass[11pt]{article}
\usepackage[plain]{fullpage}
\usepackage{amsthm, amsfonts,amsmath,amssymb}
\usepackage{xcolor}[dvipsnames]
\usepackage{tikz}
\usepackage{bbold}
\usepackage{caption}
\usepackage{hyperref}
\usepackage{thm-restate}
\usepackage{enumitem}
\usepackage{authblk}

\usetikzlibrary{arrows}

\newtheorem{theorem}{Theorem}
\newtheorem{lemma}[theorem]{Lemma}

\newtheorem{definition}[theorem]{Definition}

\newcommand{\dic}{\vec{\chi}}
\newcommand{\bid}{\overleftrightarrow}
\newcommand{\ind}[1]{[#1]}
\DeclareMathOperator{\UG}{UG}
\DeclareMathOperator\NN{\mathbb{N}}

\definecolor{g-green}{rgb}{0.235, 0.659, 0.322}
\definecolor{g-blue}{rgb}{0.0, 0.5, 1.0}

\let\leq\leqslant
\let\geq\geqslant
\let\emptyset\varnothing

\begin{document}
	
\title{Brooks-type colourings of digraphs in linear time}

\author[1]{Daniel Gonçalves}
\author[2]{Lucas Picasarri-Arrieta\footnote{Research supported by grants ANR-19-CE48-0013, ANR-17-EURE-0004, and JST ASPIRE JPMJAP2302.}$^{,}$}
\author[1]{Amadeus Reinald}

\affil[1]{LIRMM, Université de Montpellier, CNRS, Montpellier, France}
\affil[2]{National Institute of Informatics, Tokyo, Japan}

\date{}

\maketitle

\sloppy

\begin{abstract}
    
    Brooks' Theorem is a fundamental result on graph colouring, stating that the chromatic number of a graph is almost always upper bounded by its maximal degree. Lovász showed that such a colouring may then be computed in linear time when it exists. Many analogues are known for variants of (di)graph colouring, notably for list-colouring and partitions into subgraphs with prescribed degeneracy.
    One of the most general results of this kind is due to Borodin, Kostochka, and Toft, when asking for classes of colours to satisfy ``variable degeneracy" constraints. An extension of this result to digraphs has recently been proposed by Bang-Jensen, Schweser, and Stiebitz, by considering colourings as partitions into ``variable weakly degenerate" subdigraphs. Unlike earlier variants, there exists no linear-time algorithm to produce colourings for these generalisations.
    
    We introduce the notion of \emph{(variable) bidegeneracy} for digraphs, capturing multiple (di)graph degeneracy variants.
    We define the corresponding concept of $F$-dicolouring, where $F = (f_1,...,f_s)$ is a vector of functions, and an $F$-dicolouring requires vertices coloured $i$ to induce a ``strictly-$f_i$-bidegenerate" subdigraph.
    We prove an analogue of Brooks' theorem for $F$-dicolouring, generalising the result of Bang-Jensen et al., and earlier analogues in turn. 
    
    Our new approach provides a linear-time algorithm that, given a digraph $D$, either produces an $F$-dicolouring of $D$, or correctly certifies that none exist. 
    This yields the first linear-time algorithms to compute (di)colourings corresponding to the aforementioned generalisations of Brooks' theorem. In turn, it gives an unified framework to compute such colourings for various intermediate generalisations of Brooks' theorem such as list-(di)colouring and partitioning into (variable) degenerate sub(di)graphs.
\end{abstract}

\sloppy

\section{Introduction}

A $k$-colouring of a graph $G$ is a function $\alpha : V(G) \to [k]$, where for any integer $k$ we let $[k]$ be the set of integers $\{1,...,k\}$. It is proper if $\alpha^{-1}(i)$ is an independent set for any $i \leq k$.
One of the most natural ways to properly colour a simple graph $G$ is to proceed greedily, considering vertices one by one, and assigning them a colour different from their already coloured neighbours.
Then, letting $\Delta(G)$ be the maximal degree of $G$, the above yields a linear-time algorithm to properly colour $G$ with $\Delta(G)+1$ colours. In 1941, Brooks~\cite{brooksMPCPS37} showed that in most cases $\Delta(G)$ colours are enough, and gave an exact characterisation of the graphs for which $\chi(G) = \Delta +1$.
\begin{theorem}[Brooks~\cite{brooksMPCPS37}]\label{thm:brooks}
    A connected graph $G$ satisfies $\chi (G) \leq \Delta(G) + 1$ and equality holds if and only if $G$ is an odd cycle or a complete graph.
\end{theorem}
The original proof of Brooks' theorem can be adapted into a quadratic-time algorithm to construct such a $\Delta$-colouring when possible. Lovász provided a simpler proof~\cite{LOVASZ1975269}, which leads to a linear-time algorithm, as explicited by Baetz and Wood~\cite{baetz2014brooks}. An alternative linear-time algorithm was also given by Skulrattanakulchai~\cite{skulrattanakulchai2002delta}.
Going back to the less restrictive case of $(\Delta+1)$-colouring, Assadi, Chen and Khanna have recently broken the linearity barrier for several models of sublinear-time algorithms~\cite{assadiSublinear2019SODA}.

Brooks' theorem stands as one of the most fundamental results in structural graph theory. We refer the interested reader to the rencent book by Stiebitz, Schweser, and Toft~\cite{Stiebitz_Schweser_Toft_2024} for a comprehensive literature. Algorithmically, it is a canonical example of restrictions under which deciding $k$-colourability is solvable in polynomial-time.
It has since been generalised in many ways, following the introduction of new notions of colourings for (directed) graphs. 
In the same vein as Brooks' theorem, we are interested in generalisations that impose constraints on the degree of a (di)graph to guarantee its colourability. Let us also mention that such results also exist when imposing edge-(arc-)connectivity restrictions instead, see~\cite{aboulker2017coloring,stiebitz2016brooks,schweser2022coloring,aboulker2023digraph}.
Giving analogues of Brooks' theorem for more general colourings is an active line of research.
The proofs of these analogues usually imply the existence of a polynomial-time algorithm to either produce the required colouring, or decide it does not exist. Nevertheless, linear-time algorithms are known only for a few cases.
The goal of this paper is to provide such linear-time algorithms for analogues based on degree constraints. In the process, we define a more general version of colouring for digraphs, and prove the corresponding analogue of Brooks' theorem.
Before doing so, let us give a bird's eye view of the analogues of Brooks' theorem we aim to generalise, discussing their time complexity along the way.

\subsection*{Undirected analogues}

One of the first analogues considered for Brooks' theorem deals with {\it list-colourings}. 
List-colouring is a generalization of graph coloring that was introduced first
by Vizing~\cite{vizingDA29}.
Given a graph $G$, a \textit{list assignment} $L$ is a function that associates a list of ``available" colours to every vertex of $G$. An {\it $L$-colouring} of $G$ is then a proper colouring $\alpha$ of $G$ such that $\alpha(v) \in L(v)$ for every vertex $v$ of $G$.
The natural extension of Brooks' theorem in this case is to ask when $|L(v)| \geq d(v)$ for all $v$ is a sufficient condition for an $L$-colouring. Such a characterisation was given independently in 1979 by Borodin~\cite{borodinThesis}, and Erd\H{o}s, Rubin and Taylor~\cite{erdosCN26}.
The graphs for which this condition is not sufficient are \textit{Gallai trees}, which are connected graphs in which every maximal biconnected subgraph, or block, is either a complete graph or an odd cycle.
\begin{theorem}[Borodin~\cite{borodinThesis} ; Erd\H{o}s, Rubin, Taylor ~\cite{erdosCN26}]
    \label{thm:undirected_gallai}
    Let $G=(V,E)$ be a connected graph and $L$ be a list assignment of $G$ such that, for every vertex $v\in V$, $|L(v)| \geq d(v)$. If $G$ is not $L$-colourable, then $G$ is a Gallai tree and $|L(v)| = d(v)$ for every vertex $v \in V$.
\end{theorem}
This result actually generalises Brooks' theorem by setting $L(v) = [\Delta]$ for every $v \in V(G)$. 
The procedures given in the proofs leads to polynomial time algorithms, but the first linear-time algorithm finding such a colouring was given by Skulrattanakulchai~\cite{skulrattanakulchai2002delta}.

\medskip

Another generalisation of colouring for which an analogue of Brooks' theorem has been obtained relaxes the independence condition on colour classes to a degeneracy condition. Given an integer $d$, a graph $G$ is \textit{$d$-degenerate} if every non-empty subgraph of $G$ contains a vertex of degree at most $d$. Then, independent sets are exactly the sets of vertices inducing a $0$-degenerate subgraph.
Let $P=(p_1,\ldots,p_s)$ be a sequence of positive integers. A graph $G$ is \textit{$P$-colourable} if there exists an $s$-colouring $\alpha$ of $G$ such that, for every $i\in[s]$, the subgraph of $G$ induced by the colour class $\alpha^{-1}(i)$ is $(p_i-1)$-degenerate. When $p_1=\ldots=p_s=1$, observe that a $P$-colouring is exactly a proper $s$-colouring.
Bollob\'as and Manvel~\cite{bollobasBLMS11}, and Borodin~\cite{borodinDA28} independently proved the following.
\begin{theorem}[Bollob\'as and Manvel~\cite{bollobasBLMS11} ; Borodin~\cite{borodinDA28}]
    \label{thm:undirected_degeneracy}
    Let $G$ be a connected graph with maximum degree $\Delta$ and $P=(p_1, \ldots, p_s)$ be a sequence of $s\geq 2$ positive integers such that $\sum_{i=1}^s{p_i} \geq \Delta$. 
    If $G$ is not $P$-colourable, then $\sum_{i=1}^s{p_i} = \Delta$ and $G$ is a complete graph or an odd cycle. 
\end{theorem}
Brooks' theorem follows from Theorem~\ref{thm:undirected_degeneracy} by setting $s=\Delta$ and $p_i=1$ for every $i\in[\Delta]$.
In terms of complexity, the proofs of~\cite{borodinDA28} and~\cite{bollobasBLMS11} provide at best cubic algorithms, which were recently improved to a linear-time algorithm by Corsini {\it et al.}~\cite{corsiniEJC114}.

\medskip

Given these two independent generalisations of Brooks' theorem, one can naturally ask for an even more general theorem subsuming both Theorems~\ref{thm:undirected_gallai} and Theorem~\ref{thm:undirected_degeneracy}. Such a result has been obtained by Borodin, Kostochka, and Toft through the introduction of \emph{variable degeneracy}~\cite{borodinDM214}. 

Let $G=(V,E)$ be a graph and $f : V \xrightarrow[]{} \mathbb{N}$ be an integer-valued function. We say that $G$ is \textit{strictly-$f$-degenerate} if every subgraph $H$ of $G$ contains a vertex $v$ satisfying $d_H(v) < f(v)$. 
Let $s\geq 1$ be an integer and $F=(f_1,\ldots,f_s)$ be a sequence of integer-valued functions. The graph $G$ is \textit{$F$-colourable} if there exists an $s$-colouring $\alpha$ of $G$ such that, for every $i\in[s]$, the subgraph of $G$ induced by the vertices coloured $i$ is strictly-$f_i$-degenerate.
The authors consider instances $(G,F)$ where for every vertex $v \in V(G)$, $\sum_{i=1}^s f_i(v)$, which can be seen as a ``colour budget'', is at least $d(v)$. They give an explicit recursive characterisation of \emph{hard pairs}, instances satisfying the above condition that are not $F$-colourable. Informally, these instances are a generalisation of Gallai trees, where some ``monochromatic'' blocks are allowed, and in which the condition above is tight in every vertex. They then prove the following analogue, generalising both Theorems~\ref{thm:undirected_gallai} and~\ref{thm:undirected_degeneracy}.

\begin{theorem}[Borodin, Kostochka, Toft~\cite{borodinDM214}]
    \label{thm:undirected_variable_degeneracy}
    Let $G$ be a connected graph and $F=(f_1, \ldots, f_s)$ be a sequence of integer-valued functions  such that, for every vertex $v\in V(G)$, $\sum_{i=1}^sf_i(v) \geq d(v)$. Then $G$ is $F$-colourable if and only if $(G,F)$ is not a hard pair.
\end{theorem}
Note that Theorem~\ref{thm:undirected_degeneracy} can be obtained from the result above by setting $f_i$ to the constant function equal to $p_i$ for every $i\in[s]$. On the other hand, given a graph $G$ and a list assignment $L$ of $G$, one can set $f_i(v)$ to $1$ when $i\in L(v)$ and $0$ otherwise to obtain Theorem~\ref{thm:undirected_gallai}.
We mention that Theorem~\ref{thm:undirected_variable_degeneracy} has been extended to hypergraphs, and refer the interested reader to~\cite{schweserJGT96,schweserEJGTA9}. It has also been generalised to correspondence colouring (which is also known as DP-colouring), see~\cite{kostochkaDM346}.

Recall linear-time algorithms have been obtained to construct the colourings guaranteed by Theorem~\ref{thm:undirected_degeneracy} and Theorem~\ref{thm:undirected_gallai}. Nevertheless, to the best of the authors' knowledge, until this paper, no linear-time algorithm existed for constructing the colourings of Theorem~\ref{thm:undirected_variable_degeneracy}, which was only shown to be polynomial.

\subsection*{Directed analogues}

In 1982, Neumann-Lara~\cite{neumannlaraJCT33} introduced the notions of dicolouring for digraphs, generalising proper colourings of graphs. 
A $k$-colouring of a digraph $D$ is a \textit{$k$-dicolouring} if vertices coloured $i$ induce an acyclic subdigraph of $D$ for each $i \in [k]$. The \textit{dichromatic number} of $D$, denoted by $\dic(D)$, is the smallest $k$ such that $D$ admits a $k$-dicolouring. 
There is a one-to-one correspondence between the proper $k$-colourings of a graph $G$ and the $k$-dicolourings of its associated bidirected graph $\bid{G}$ (where all edges of $G$ are replaced by digons), and in particular $\chi(G) = \dic(\bid{G})$. Hence, every result on graph proper colourings can be seen as a result on dicolourings of bidirected graphs, and it is natural to study whether the result can be extended to all digraphs. 
The following directed version of Brooks' theorem was first obtained by Mohar~\cite{moharLAA432}. Aboulker and Aubian~\cite{aboulkerDM113193} recently gave four new proofs of the result.
\begin{theorem}[Mohar~\cite{moharLAA432}]
    \label{thm:brooks_directed}
    Let $D=(V,A)$ be a connected digraph. Then 
    \[ \dic(D) \leq \max \{ d^+(v),d^-(v) \mid v\in V\} +1\]
    and equality holds if and only if $D$ is a directed cycle, a bidirected odd cycle, or a bidirected complete graph.
\end{theorem}
The constructive proofs of the result yields a linear-time algorithm, and it may be possible to derive a linear-time algorithm from the generalisation of Lovász's proof given in~\cite{aboulkerDM113193}.
In the special case of oriented graphs (\textit{i.e.} digraphs of girth at least three), the second author~\cite{picasarriJGT106} strengthened Theorem~\ref{thm:brooks_directed} and proved that $\dic(D) \leq \Delta_{\min} = \max_v \min(d^-(v),d^+(v))$ holds unless $\Delta_{\min} \leq 1$. The deterministic proof gives a linear-time algorithm for finding a dicolouring of $D$ using at most $\Delta_{\min}$ colours.

Harutyunyan and Mohar~\cite{harutyunyanSIDMA25} proved an analogue for list-dicolouring of digraphs.  
A \textit{directed Gallai tree} is a digraph in which every block is a directed cycle, a bidirected odd cycle or a bidirected complete graph. Given a list assignment of a digraph $D$, an $L$-dicolouring $\alpha$ is a dicolouring of $D$ such that $\alpha(v) \in L(v)$ holds for every vertex $v$ of $D$. 
\begin{theorem}[Harutyunyan and Mohar~\cite{harutyunyanSIDMA25}]\label{thm:gallai_directed}
    Let $D=(V,A)$ be a connected digraph and $L$ a list assignment of $D$ such that $|L(v)| \geq \max(d^+(v),d^-(v))$ holds for every vertex $v\in V$. If $D$ is not $L$-dicolourable, then $D$ is a directed Gallai tree.
\end{theorem}
This also generalises Theorem~\ref{thm:undirected_gallai}, and it is shown in~\cite{harutyunyanSIDMA25} that such a colouring may be found in linear time when it exists. 

In the current direction, the most general analogue of Brooks' theorem was introduced in 2020 by Bang-Jensen, Schweser, and Stiebitz~\cite{bangSIDMA36}.
The authors extend the notion of variable degeneracy to digraphs, and define the corresponding colouring notion, dubbed \emph{$f$-partition}. Given a function $h : V(D) \xrightarrow[]{} \mathbb{N}$, a digraph $D$ is $h$-min-degenerate\footnote{Actually, in~\cite{bokalWeakDegeneracy,bangSIDMA36} (variable) min-degeneracy is called \emph{(variable) weak degeneracy}, but these names are better suited for our presentation.} if for every non-empty subdigraph of $D$, there exists some $v \in V(D)$ such that $\min(d^-(v),d^+(v)) < h(v)$.
Then, letting $f = (f_1,...,f_s)$, where $f_i : V(D) \xrightarrow[]{} \mathbb{N}$ for any $i \in [s]$, an $f$-partition of $D$ is a colouring $\alpha$ of $V(D)$ such that vertices coloured $i$ induce an $f_i$-min-degenerate subgraph for every $i$.
Similar to the case of variable degeneracy, the instances considered are of the form $(D,f)$ such that $\sum_{i=1}^s f_i(v) \geq \max(d^-(v),d^+(v))$ for every $v \in V(D)$. Then, the authors define hard pairs in an analogue manner to~\cite{borodinDM214}, and show those are exactly the instances that are not $f$-partionable.
\begin{theorem}[Bang-Jensen, Schweser, Stiebitz~\cite{bangSIDMA36}]
    \label{thm:f-partition}
    Let $D$ be a connected digraph, let $s \geq 1$ be an integer, and let $f : V (D) \xrightarrow[]{} \mathbb{N^s}$, be a vector function such that $\sum_{i=1}^s f_i(v) \geq \max(d^-(v),d^+(v))$ for all $v \in V(D)$. Then, $D$ is $f$-partitionable if and only if $(D,f)$ is not a hard pair.
\end{theorem}
This notion generalises both dicolouring and list-dicolouring, as well as all variants of undirected colouring presented above. In particular, all the analogues of Brooks' theorem presented until now can be seen as a particular case of Theorem~\ref{thm:f-partition}.
The authors provide algorithms to produce such colourings or decide of a hard pair, and argue they are polynomial. The time-complexity is not given explicitly, but appears to be at least quadratic.

\subsection*{Our contributions}

Our goal is to obtain a linear-time algorithm to construct the colourings guaranteed by the analogues of Brooks' theorem introduced earlier, or decide they are not colourable. The cases for which such an algorithm was missing being Theorem~\ref{thm:undirected_variable_degeneracy} and Theorem~\ref{thm:f-partition}.
As the latter is the most general, it is natural to look for a linear-time algorithm to find the $f$-partitions considered in~\cite{bangSIDMA36}.
Towards this, we define the more expressive notion of variable bidegeneracy for digraphs, as well as the corresponding $F$-dicolouring concept. We obtain an analogue of Brooks' theorem for $F$-dicolouring, as well as a linear-time algorithm that either produces such a colouring, or certifies that there is none. Our strategy does not follow previous approaches in the literature, and the novel point of view given by $F$-dicolourings is what ultimately leads to linear-time algorithms for earlier notions as well.

\subsubsection*{Variable bidegeneracy}

The key to our generalisation of Brooks' theorem and corresponding colouring algorithm is the new concept of (variable) \emph{bidegeneracy} for digraphs.
Defining degeneracy notions for digraphs can be done in various ways by imposing conditions on the in- and out-degrees. The first such example in the literature is given by Bokal {\it et al.}~\cite{bokalWeakDegeneracy}, who introduced min-degeneracy. Recall that a digraph is $k$-min-degenerate if all its subdigraphs admit a vertex such that the minimum of the in- and out-degrees is at most $k$. One may also ask for only the in-degree (or only the out-degree) to be at most $k$, which corresponds to the $k$-in-degeneracy (or $k$-out-degeneracy) of Bousquet {\it et al.}~\cite{bousquetEJC116}.
We introduce the concept of bidegeneracy for digraphs. For any pair of integers $(p,q)$, a digraph $D$ is {\it strictly-$(p,q)$-bidegenerate} if for every subdigraph of $D$, there exists a vertex $v$ such that either $d^-(v) < p$ or $d^+(v) < q$.

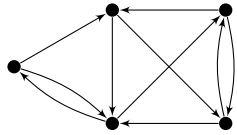
\begin{figure}[hbtp]
    \begin{minipage}{\linewidth}
        \begin{center}	
            \begin{tikzpicture}[scale=1, every node/.style={transform shape}]
                \tikzset{vertex/.style = {circle,fill=black,minimum size=5pt, inner sep=0pt}}
                \tikzset{edge/.style = {->,> = latex'}}
                
                \node[vertex] (1) at (-1.3,0.75) {};
                \node[vertex] (2) at (0,1.5) {};
                \node[vertex] (3) at (0,0) {};
                \node[vertex] (4) at (1.5,0) {};
                \node[vertex] (5) at (1.5,1.5) {};
                
                \draw[edge] (1) to (2) {};
                \draw[edge, bend left=13] (1) to (3) {};
                \draw[edge] (2) to (3) {};
                \draw[edge,] (2) to (4) {};
                \draw[edge] (3) to (5) {};
                \draw[edge, bend left=13] (3) to (1) {};
                \draw[edge] (4) to (3) {};
                \draw[edge, bend left=13] (4) to (5) {};
                \draw[edge, bend left=13] (5) to (4) {};
                \draw[edge] (5) to (2) {};
            \end{tikzpicture}
            \caption{A digraph which is strictly-$(2,0)$-bidegenerate while it is not strictly-$(0,2)$-bidegenerate.}
            \label{fig:(2_0)-degenerate}
        \end{center}    
    \end{minipage}
\end{figure}

Bidegeneracy allows us to express both min-degeneracy as well as in- or out-degeneracy.
Indeed, observe that a digraph $D$ is $k$-min-degenerate for some $k \in \mathbb{N}$ if and only if it is strictly-$(k+1,k+1)$-bidegenerate.
Then, $D$ is $k$-in-degenerate (respectively $k$-out-degenerate) if and only if it is strictly-$(k+1,0)$-bidegenerate (respectively strictly-$(0,k+1)$-bidegenerate).
In particular, for a digraph, it is equivalent to be acyclic, to be $0$-min-degenerate, or to be strictly-$(1,1)$-bidegenerate.
Furthermore, bidegeneracy constraints express a finer range of classes of digraphs, lying strictly between classes defined by other notions.
For instance, the digraph illustrated in Figure~\ref{fig:(2_0)-degenerate} is strictly-$(2,0)$-bidegenerate, thus $2$-min-degenerate, yet it is not strictly-$(0,2)$-bidegenerate.
Indeed, bidegeneracy is less symmetric in nature, and importantly, it is not preserved under reversal of all arcs in a given digraph, while min-degeneracy is.

Given any degeneracy notion for (di)graphs, the corresponding variable notion may be defined through a function on the vertices, imposing different restrictions on each vertex of the (di)graph at hand.
Notably, the degeneracy of~\cite{bangSIDMA36} corresponds to the variable version of the degeneracy introduced in~\cite{bokalWeakDegeneracy}. To the best of our knowledge, there are no known results on the variable in-degeneracy (or variable out-degeneracy).
We now define variable bidegeneracy, for a digraph $D$, through a function $f : V(D) \xrightarrow[]{} \mathbb{N}^2$. We denote the projection of $f$ on the first and second coordinates by $f^-$ and $f^+$ respectively. Then, $D$ is said to be \emph{strictly-$f$-bidegenerate} if for every subdigraph $H$ of $D$, there exists some $v \in V(D)$ such that either $d_H^-(v) < f^-(v)$ or $d_H^+(v) < f^+(v)$.
If $f^-(v)=f^+(v)$ for every $v\in V(D)$, $f$ is said to be symmetric. 
Again, for a function $h$, a digraph is $h$-min-degenerate if and only if it is strictly-$(h,h)$-bidegenerate.
As a consequence, our notion also generalises all the degeneracy notions used for the analogues of Brooks' theorem presented here.

\subsubsection*{$F$-dicolouring}

With this in hand, we are ready to define our corresponding notion of dicolourings.
Let $D$ be a digraph, $s$ be a positive integer, and $F = (f_1,\ldots,f_s)$ be a sequence of functions $f_i : V(D) \xrightarrow{} \NN^2$. Then, $D$ is \textit{$F$-dicolourable} if there exists an $s$-colouring $\alpha$ of $V(D)$ such that, for every $i\in[s]$, the subdigraph of $D$ induced by the vertices coloured $i$ is strictly-$f_i$-bidegenerate. Such a colouring $\alpha$ is called an \textit{$F$-dicolouring}.

Deciding if a digraph is $F$-dicolourable is clearly NP-hard because it includes a large collection of NP-hard problems for specific values of $F$. For instance, deciding if a digraph $D$ has dichromatic number at most two, which is shown to be NP-hard in~\cite{chenSIC37}, consists exactly in deciding whether $D$ is $((1,0),(1,0))$-dicolourable (when considering $(1,0)$ as a constant function).
We thus restrict ourselves to \emph{valid pairs} $(D,F)$, for which $D$ is connected and the following holds:
\begin{equation}\label{property_D_F}
    \forall v\in V(D), \text{~~~} \sum_{i=1}^s f_i^-(v) \geq d^-(v) \mbox{~~~and~~~} \sum_{i=1}^s f_i^+(v) \geq d^+(v).\tag{$\star$}
\end{equation}
This is a natural condition on the ``colour budget'' $F$ in terms of the degrees of $D$, and directly generalises analogues presented earlier. Lowering this budget by only one would allow the expression of hard instances, such as 3-colourability of graphs with maximum degree 4, which is already NP-hard (see~\cite{garey1979}).

In Subsection~\ref{subsec:hard-pairs}, we define hard pairs $(D,F)$, consisting of a digraph $D$ and a sequence $F$ of functions satisfying Property~\eqref{property_D_F}. The definition is inductive, with the base cases consisting of hard blocks (for which $D$ is biconnected), which may then be glued together by identifying vertices to form other hard pairs. These turn out to be a straightforward adaptation of the hard pairs of~\cite{bangSIDMA36}, and capture other notions of ``hard'' instances (such as Gallai trees) for the previous analogues of Brooks' theorem.
We generalise Brooks' theorem for valid pairs, by showing they are $F$-dicolourable if and only if they are not hard.
\begin{restatable}{theorem}{bivariable}
    \label{thm:bivariable}
    Let $(D,F)$ be a valid pair. Then $D$ is $F$-dicolourable if and only if $(D,F)$ is not a hard pair (as defined in Subsection~\ref{subsec:hard-pairs}). Moreover, there is an algorithm running in time $O(|V(D)|+|A(D)|)$ that decides if $(D,F)$ is a hard pair, and that outputs an $F$-dicolouring if it is not.
\end{restatable}
While the theorem is of independent interest for the new notion of variable bidegeneracy, it is an unifying result allowing for linear-time algorithms for all the analogues above. Importantly, it yields the first linear-time algorithm for the $F$-colourability of Theorem~\ref{thm:undirected_variable_degeneracy} and the $f$-partitioning of Theorem~\ref{thm:f-partition}. 
Note also that our complexity does not depend on the total number of colours $s$, such as the algorithm of~\cite{skulrattanakulchai2002delta} for Theorem~\ref{thm:undirected_gallai}. As the input has size $O(|V(D)|+|A(D)|+s\cdot |V(D)|)$, our algorithm can thus be sublinear in the input size if $(|V(D)|+|A(D)|)$ is asymptotically dominated by $(s\cdot |V(D)|)$.
From now on, ``linear-time complexity'' means linear in the number of vertices and arcs of the considered digraph.

We begin with a few preliminaries and the definition of hard pairs in Section~\ref{section:brooks_digraphs}. Then, Section~\ref{section:proof_thm_bivariable} is devoted to the proof of Theorem~\ref{thm:bivariable}.

\section{Preliminaries}
\label{section:brooks_digraphs}

We refer the reader to~\cite{bang2009} for notation and terminology on digraphs not explicitly defined in this paper.
Let $D=(V,A)$ be a digraph. The vertex set of $D$ and its arc set are denoted respectively by $V(D)$ and $A(D)$. A \textit{digon} is a pair of arcs in opposite directions between the same vertices. A \textit{simple arc} is an arc which is not in a digon. 
The \textit{bidirected graph} associated with a graph $G$, denoted by $\bid{G}$, is the digraph obtained from $G$ by replacing every edge by a digon. The \textit{underlying graph} of $D$, denoted by $\UG(D)$, is the undirected graph $G$ with vertex set $V(D)$ in which $uv$ is an edge if and only if $uv$ or $vu$ is an arc of $D$. 
A digraph is {\it connected} if its underlying graph is connected.
A connected (di)graph $D$ is \textit{biconnected} if it has at least two vertices and it remains connected after the removal of any vertex. Note that $K_2$ is biconnected. A \textit{block} of $D$ is a maximal biconnected sub(di)graph of $D$. A {\it cut-vertex} of $D$ is a vertex $x\in V(D)$ such that $D-x$ is disconnected. An \textit{end-block} is a block with at most one cut-vertex.
Let $v$ be a vertex of a digraph $D$. The {\it out-degree} (respectively {\it in-degree}) of $v$, denoted by $d^+(v)$ (respectively $d^-(v)$), is the number of arcs leaving (resp. entering) $v$. The \textit{degree} of $v$, denoted by $d(v)$, is the sum of its in- and out-degrees. A \textit{directed cycle} is a connected digraph in which every vertex $v$ satisfies $d^+(v) = d^-(v) = 1$. 
A digraph is \textit{acyclic} if it does not contain any directed cycle.
Recall, a colouring $\alpha$ of the vertices of a digraph $G$ is considered as a function $\alpha : V(G) \xrightarrow[]{} [k]$, where $[k]$ denotes the set of integers $\{1,..,k\}$.

For two elements $p=(p_1,p_2)$ $q=(q_1,q_2)$ of $\NN^2$, we denote by $p \leq q$ the relation $p_1 \leq q_1$ and $p_2 \leq q_2$. We denote by $p<q$ the relation $p\leq q$ and $p\neq q$. 

\subsection{Hard Pairs}\label{subsec:hard-pairs}

We are now ready to define \textit{hard pairs}, which are pairs $(D,F)$ that satisfy Property~\eqref{property_D_F} tightly, with specific conditions on $F$, and which admit a certain block tree structure on $D$.
We say that $(D,F)$ is a \textit{hard pair} if one of the following four conditions holds:

\begin{enumerate}[label=$(\roman*)$]
    \item $D$ is a biconnected digraph and there exists $i\in [s]$ such that, for every vertex $v\in V$, $f_i(v)=(d^-(v),d^+(v))$ and $f_k(v) = (0,0)$ when $k\neq i$. 
    
    We refer to such a hard pair as a \textit{monochromatic hard pair}.
    
    \item $D$ is a bidirected odd cycle and the functions $f_1,\ldots,f_s$ are all constant equal to $(0,0)$ except exactly two that are constant equal to $(1,1)$.
    
    We refer to such a hard pair as a \textit{bicycle hard pair}.
    
    \item $D$ is a bidirected complete graph, the functions $f_1,\ldots,f_s$ are all constant and symmetric, and for every vertex $v$ we have  $\sum_{i=1}^sf_i^+(v) = |V(D)|-1$.
    
    We refer to such a hard pair as a \textit{complete hard pair}.
    
    \item $(D,F)$ is obtained from two hard pairs $(D^1,F^1)$ and $(D^2,F^2)$ by identifying two vertices $x_1 \in V(D^1)$ and $x_2\in V(D^2)$ into a new vertex $x\in V(D)$, such that for every vertex $v\in V(D)$ we have: 
    \[
    \mbox{~~~~~~~}f_k(v) = \left\{
    \begin{array}{lll}
        f^1_k(v) & \mbox{if } v\in V(D^1) \setminus \{x_1\} \\
        f^2_k(v) & \mbox{if } v\in V(D^2) \setminus \{x_2\} \\
        f^1_k(v) + f^2_k(v) & \mbox{if } v = x.
    \end{array}
    \right.
    \]
    where $F^1 = (f^1_1,\ldots,f^1_s)$ and $F^2 = (f^2_1,\ldots,f^2_s)$.
    
    We refer to such a hard pair as a \textit{join hard pair}.
\end{enumerate}

\begin{figure}[hbtp]
    \begin{minipage}{\linewidth}
        \begin{center}	
            \begin{tikzpicture}[scale=1, every node/.style={transform shape}]
                \tikzset{vertex/.style = {circle,fill=black,minimum size=5pt, inner sep=0pt}}
                \tikzset{edge/.style = {->,> = latex}}
                
                \node[vertex, label=left:${\textcolor{purple}{(1,1)}, \textcolor{g-green}{(2,2)}}$] (v1) at  (-10*9/10,0*9/10) {};
                \node[vertex, label=above:${\textcolor{purple}{(1,1)},\textcolor{g-green}{(2,2)}}$] (v2) at  (-8*9/10,2*9/10) {};
                \node[vertex, label=below:${\textcolor{purple}{(1,1)},\textcolor{g-green}{(2,2)}}$] (v3) at  (-8*9/10,-2*9/10) {};
                \node[vertex, label=right:${\textcolor{purple}{(1,3)},\textcolor{g-green}{(2,2)}}$] (v4) at  (-6*9/10,0*9/10) {};
                \foreach \i in {1,...,3}
                {{\pgfmathtruncatemacro{\j}{\i + 1}
                        \foreach \k in {\j,...,4}{   
                            \draw[edge, bend left=10] (v\i) to (v\k);
                            \draw[edge, bend left=10] (v\k) to (v\i);
                }}}
                \node[vertex, label=above:${\textcolor{purple}{(2,1)}}$] (v5) at  (-4*9/10,2*9/10) {};
                \node[vertex, label=below:${\textcolor{purple}{(1,2)}}$] (v6) at  (-4*9/10,-2*9/10) {};
                \draw[edge] (v6) to (v5);
                \draw[edge] (v4) to (v5);
                \draw[edge] (v4) to (v6);
                
                \node[vertex, label=right:${\textcolor{purple}{(3,1)},\textcolor{g-blue}{(1,1)}}$] (v7) at (180:2*9/10){};
                \node[vertex, label=below:${\textcolor{purple}{(1,1)},\textcolor{g-blue}{(1,1)}}$] (v8) at (-108:2*9/10){};
                \node[vertex, label=right:${\textcolor{purple}{(1,1)},\textcolor{g-blue}{(1,1)}}$] (v9) at (-36:2*9/10){};
                \node[vertex, label=right:${\textcolor{purple}{(1,1)},\textcolor{g-blue}{(1,1)}}$] (v10) at (36:2*9/10){};
                \node[vertex, label=above:${\textcolor{purple}{(1,1)},\textcolor{g-blue}{(1,1)}}$] (v11) at (108:2*9/10){};
                \foreach \i in {7,...,10}
                {{\pgfmathtruncatemacro{\j}{\i + 1}
                        \draw[edge, bend left=10] (v\i) to (v\j);
                        \draw[edge, bend left=10] (v\j) to (v\i);
                }}
                \draw[edge, bend left=10] (v11) to (v7);
                \draw[edge, bend left=10] (v7) to (v11);
                \draw[edge] (v6) to (v7);
                \draw[edge] (v5) to (v7);
            \end{tikzpicture}
            \caption{An example of a hard pair $(D,F)$ where $F$ is a sequence of three functions $f_1,f_2,f_3$ the values of which are respectively represented in red, green, and blue. For the sake of readability, the value of $f_i$ is missing when it is equal to $(0,0)$.}
            \label{fig:hard-pair}
        \end{center}    
    \end{minipage}
\end{figure}
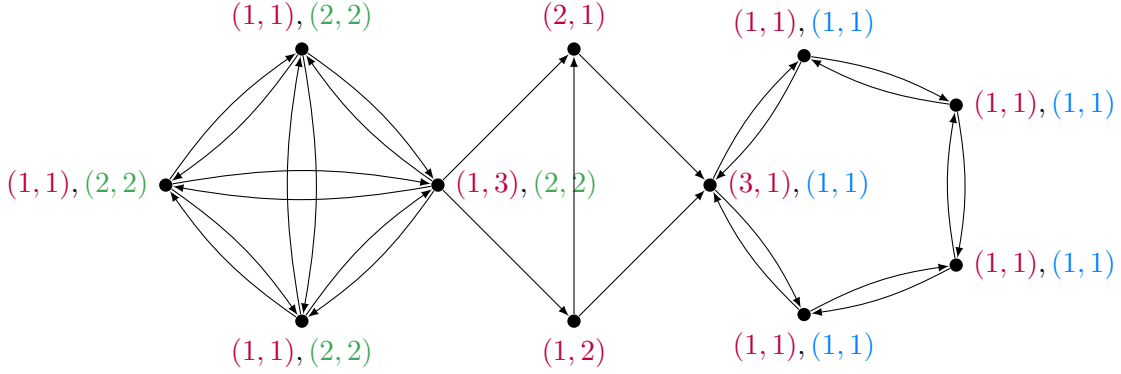

See Figure~\ref{fig:hard-pair} for an illustration of a hard pair. 

In the following, we give a short proof that if $(D,F)$ is a hard pair then $D$ is not $F$-dicolourable. 

\subsection{Hard pairs are not dicolourable}
\label{section:hard_pair_not_dicolourable}

\begin{lemma}\label{lemma:hard_pairs_not_dicolourable}
    Let $(D, F)$ be a hard pair, then $D$ is not $F$-dicolourable.
\end{lemma}
\begin{proof}
    We proceed by induction.
    Let $(D,F=(f_1,\ldots,f_s))$ be a hard pair. Assume for a contradiction that it admits an $F$-dicolouring $\alpha$ with colour classes $V_1,\ldots,V_s$. We distinguish four cases, depending on the kind of hard pair $(D,F)$ is.
    \begin{enumerate}[label=$(\roman*)$]
        \item If $(D,F)$ is a monochromatic hard pair, let $i\in [s]$ be such that, for every vertex $v\in V$, $f_i(v) = (d^-(v),d^+(v))$ and $f_k(v) = (0,0)$ when $k\neq i$. Then, for every $k\neq i$, $V_k$ must be empty, for otherwise $D\ind{V_k}$ is not strictly-$f_k$-bidegenerate.
        Therefore, we have $V=V_i$ and $D$ must be strictly-$f_i$-bidegenerate. Hence $D$ contains a vertex $v$ such that $d^-(v) < f_i^-(v)$ or $d^+(v) < f_i^+(v)$, a contradiction.
        
        \item If $(D,F)$ is a bicycle hard pair, let $i,j\in [s]$ be distinct integers such that, for every vertex $v\in V$, $f_i(v) = f_j(v) = (1,1)$ and $f_k(v) = (0,0)$ when $k\notin \{i,j\}$. Again, $V_k$ must be empty when $k\notin \{i,j\}$, so $(V_i,V_j)$ partitions $V$. Since $D$ is a bidirected odd cycle, it is not bipartite, so $D\ind{V_i}$ or $D\ind{V_j}$ contains a digon between vertices $u$ and $v$. Assume that $D\ind{V_i}$ does, then $H=D\ind{\{u,v\}}$ must contain a vertex $x$ such that $d^-_H(x) < f_i^-(x)$ or $d^+_H(x) < f_i^+(x)$, a contradiction since $f_i(x) = (d^-_H(x),d^+_H(x)) = (1,1)$ for every $x\in \{u,v\}$.
        
        \item If $(D,F)$ is a complete hard pair, then for each $i\in [s]$, since $D\ind{V_i}$ is strictly-$f_i$-bidegenerate, we have $|V_i| \leq f_i^-(v)$ where $v$ is any vertex (recall that the functions $f_1,\ldots,f_s$ are constant and symmetric). Then $\sum_{i=1}^s|V_i|\leq \sum_{i=1}^sf_i^-(v) = |V(D)|-1$, a contradiction since $(V_1,\ldots,V_s)$ partitions $V$.
        
        \item Finally, if $(D,F)$ is a join hard pair,
        $(D,F)$ is obtained from two hard pairs $(D^1,F^1)$ and $(D^2,F^2)$ by identifying two vertices $x_1 \in V(D^1)$ and $x_2\in V(D^2)$ into a new vertex $x\in V(D)$, such that for every vertex $v\in V(D)$ we have: 
        \[
        f_k(v) = \left\{
        \begin{array}{lll}
            f^1_k(v) & \mbox{if } v\in V(D^1) \setminus \{x_1\} \\
            f^2_k(v) & \mbox{if } v\in V(D^2) \setminus \{x_2\} \\
            f^1_k(v) + f^2_k(v) & \mbox{if } v = x.
        \end{array}
        \right.
        \]
        where $F^1 = (f^1_1,\ldots,f^1_s)$ and $F^2 = (f^2_1,\ldots,f^2_s)$.
        
        By induction, we may assume that $D^1$ is not $F^1$-dicolourable and $D^2$ is not $F^2$-dicolourable. 
        Let $\alpha^1$ and $\alpha^2$ be the following $s$-colourings of $D^1$ and $D^2$.
        \[
        \alpha^1(v) = \left\{
        \begin{array}{lll}
            \alpha(v) \text{ if $v \neq x_1$} \\
            \alpha(x) \text{ otherwise}.
        \end{array}
        \right.
        ~~~\text{and}~~~
        \alpha^2(v) = \left\{
        \begin{array}{lll}
            \alpha(v) \text{ if $v \neq x_2$} \\
            \alpha(x) \text{ otherwise}.
        \end{array}
        \right.
        \]
        
        We denote respectively by $U_1,\ldots,U_s$ and $W_1,\ldots,W_s$ the colour classes of $\alpha^1$ and $\alpha^2$.
        Since $D^1$ is not $F^1$-dicolourable, there exist $i\in [s]$ and a subdigraph $H^1$ of $D^1\ind{U_i}$ such that every vertex $u\in V(H^1)$ satisfies $d^-_{H^1}(u) \geq {f^{1}_i}^-(u)$ and $d^+_{H^1}(u) \geq {f^{1}_i}^+(u)$.
        We claim that $x_1\in V(H^1)$. Assume not, then $H^1$ is a subdigraph of $D\ind{V_i}$ and every vertex $u$ in $H^1$ satisfies $f_i(u) = f_i^1(u)$, and so it satisfies $d^-_{H^1}(u) \geq f_i^-(u)$ and $d^+_{H^1}(u) \geq f_i^+(u)$. This is a contradiction to $D\ind{V_i}$ being strictly-$f_i$-bidegenerate.
        
        Hence, we know that $x_1 \in V(H^1)$. Swapping the roles of $D^1$ and $D^2$, there exists an index $j\in [s]$ and a subdigraph $H^2$ of $D^2\ind{W_j}$ such that every vertex $w\in V(H^2)$ satisfies $d^-_{H^2}(w) \geq {f^{2}_j}^-(w)$ and $d^+_{H^2}(w) \geq {f^{1}_j}^+(w)$. Analogously, we must have $x_2\in V(H^2)$.
        
        By definition of $\alpha^1$ and $\alpha^2$, we have  $i=j$ and $x\in V_i$. Let $H$ be the subdigraph of $D\ind{V_i}$ obtained from $H_1$ and $H_2$ by identifying $x_1$ and $x_2$ into $x$.
        For every vertex $u\in V(H) \cap V(H^1)$, we have $d_H^-(u) = d^-_{H^1}(u) \geq f_i^-(u)$ and $d_H^+(u) = d^+_{H^1}(u) \geq f_i^+(u)$. Analogously, for every vertex $w\in V(H) \cap V(H^2)$, we have $d^-_H(w) = d^-_{H^2}(w) \geq {f^{2}_i}^-(w)$ and $d^+_H(w) = d^+_{H^2}(w) \geq {f^{1}_i}^+(w)$. 
        Finally, we have 
        \begin{align*}
            d^-_H(x) &= d^-_{H^1}(x_1) + d^-_{H^2}(x_2) \geq {f^1_i}^-(x_1) +{f^2_i}^-(x_2) = f_i^-(x)\\
            \text{and~~~}d^+_H(x) &= d^+_{H^1}(x_1) + d^+_{H^2}(x_2) \geq {f^1_i}^+(x_1) +{f^2_i}^+(x_2) = f_i^+(x).
        \end{align*}
        This is a contradiction to $D\ind{V_i}$ being strictly-$f_i$-bidegenerate.
        \qedhere
    \end{enumerate}
\end{proof}

\section{Dicolouring non-hard pairs in linear time}
\label{section:proof_thm_bivariable}

In this section, we prove Theorem~\ref{thm:bivariable} by describing the corresponding algorithm, proving its correctness, and its time complexity. 
Our algorithm starts by testing whether $(D,F)$ is a hard pair. If it is not, we pipeline various algorithms reducing the initial instance, while partially colouring the vertices pruned along the way, until a solution is found.

We first discuss the data structure encoding the input in Subsection~\ref{subsec:data_structure}.
We give some technical definitions used in the proof and some preliminary remarks in Subsection~\ref{subsec:preliminaries}.
In the following subsections we describe different reduction steps unfolding in our algorithm. We consider valid pairs for which the constraints are ``loose'' in Subsection~\ref{subsec:reduction-tightening}, based on a simple greedy algorithm presented in Subsection~\ref{subsec:greedy-algorithm}. Then, we consider ``tight'' valid pairs. In Subsection~\ref{subsec:reduction-block}, we show how to reduce a tight and valid pair $(D,F)$ into another one $(D',F')$ for which $D'$ is biconnected and smaller than $D$.
At this point, we may test the hardness of $(D',F')$ to decide whether $(D,F)$ is hard. Otherwise, the remainder of the section consists in exhibiting an $F'$-dicolouring of $D'$, yielding an $F$-dicolouring of $D$.
In Subsection~\ref{subsec:reduction-2colblocks}, we reduce the instance to one involving at most two colours, and which is still biconnected.
Then, we show how to solve biconnected instances with two colours in Subsection~\ref{subsec:algo-block-2col-particular} and in Subsection~\ref{subsec:algo-block-2col-general}, through the use of ear-decompositions. We conclude with the proof of Theorem~\ref{thm:bivariable} in Subsection~\ref{subsec:proof_thm_bivariable}.

\subsection{Data structures}\label{subsec:data_structure}

We need appropriate data structures to process the entry pair $(D,F)$.
The following structures are standard, but let us list their properties.
The digraph $D$ is encoded in space $O(|V(D)|+|A(D)|)$ with a data structure allowing, for every vertex $v$,
\begin{itemize}
    \item to access the values $|V(D)|$, $|A(D)|$, $d^-(v)$, $d^+(v)$, and $|N^-(v)\cap N^+(v)|$ in $O(1)$ time,
    \item to enumerate the vertices of the sets $V(D)$, $N^-(v)$, $N^+(v)$ in $O(1)$ time per vertex,
    \item to delete a vertex $v$ (and update all the related values and sets) in $O(d^-(v)+d^+(v))$ time, and
    \item to compute a spanning tree rooted at a specified root in $O(|V(D)|+|A(D)|)$ time.
\end{itemize}

The functions $F=(f_i)_{i\in[s]}$ are encoded in $O(s\cdot |V(D)|)$ space in a data structure allowing, for every vertex $v$,
\begin{itemize}
    \item to read or modify $f_i^-(v)$ and $f_i^+(v)$ in $O(1)$ time. 
    \item to enumerate (only) the colours $i$ such that $f_i(v) \neq (0,0)$, in $O(1)$ time per such colour. 
\end{itemize}
This data structure is simply a $s \times |V(D)|$ table, dynamically maintained, with pointers linking the cells $(i,v)$ and $(j,v)$ if $f_i(v)\neq (0,0)$, $f_j(v)\neq (0,0)$, and $f_k(v)= (0,0)$ for every $k$ such that $i<k<j$.

The output is a vertex colouring which we build along the different steps of our algorithm, and is simply encoded in a table.

\subsection{Preliminaries}\label{subsec:preliminaries}

Let $(D,F=(f_1,\ldots,f_s))$ be a (non-necessarily valid) pair, let $X\subseteq V(D)$ be a subset of vertices of $D$, and $\alpha : X \to [s]$ be a partial $s$-colouring of $D$. 
Let $D' = D-X$ and $F'=(f_1',\ldots,f_s')$ be defined as follows:
\begin{align*}
    f'^-_i(u) &= \max(0, f_i^-(u) - |\alpha^{-1}(i) \cap N^-(u)|)\\
    \text{and~~~} f'^+_i(u) &= \max(0, f_i^+(u) - |\alpha^{-1}(i) \cap N^+(u)|)
\end{align*}
We call $(D',F')$ the pair \textit{reduced} from $(D,F)$ by $\alpha$. 

\begin{lemma}\label{lemma:safe-colouring}
    Consider a (non-necessarily valid) pair $(D,F)$ and an $F$-dicolouring $\alpha$ of $D[X]$, for a subset $X\subseteq V(D)$.
    Let $(D',F')$ be the pair reduced from $(D,F)$ by $\alpha$.
    Then, combining $\alpha$ with any $F'$-dicolouring of $D'$ yields an $F$-dicolouring of $D$.
    
    Furthermore, there is an algorithm that given a pair $(D,F)$, and a colouring $\alpha$ of some set $X\subseteq V(D)$, outputs the reduced pair $(D',F')$ in $O(\sum_{v\in X} d(v))$ time.
\end{lemma}

\begin{proof}
    Let $\beta$ be any $F'$-dicolouring of $D'$. We will show that the combination $\gamma$ of $\alpha$ and $\beta$ is necessarily an $F$-dicolouring of $D$. We formally have
    \[ \gamma(v) = \left\{
    \begin{array}{lll}
        \alpha(v) & \mbox{if } v\in X \\
        \beta(v) & \mbox{otherwise.}
    \end{array}
    \right.\]
    
    Let $i\in [s]$ be any colour, we will show that $D\ind{V_i}$ is strictly-$f_i$-bidegenerate, where $V_i = \gamma^{-1}(i)$, which implies the first part of the statement. To this purpose, let $H$ be any subdigraph of $D\ind{V_i}$, we will show that $H$ contains a vertex $v$ satisfying $d^-_H(v)<f_i^-(v)$ or $ d^+_H(v) < f_i^+(v)$. 
    Observe first that, if $V(H) \subseteq X$, the existence of $v$ is guaranteed since $\alpha$ is an $F$-dicolouring of $D\ind{X}$. 
    Henceforth assume that $V(H) \setminus X \neq \emptyset$, and let $H'$ be $H - X$. Since $V(H') \subseteq \beta^{-1}(i)$, by hypothesis on $\beta$, it is strictly-$f'_i$-bidegenerate. Hence there must be a vertex $v\in V(H')$ such that $d^-_{H'}(v)<{f_i'}^-(v)$ or $ d^+_{H'}(v) < {f_i'}^+(v)$. 
    If $d^-_{H'}(v)<{f_i'}^-(v)$, we obtain
    \[ d^-_H(v) \leq d^-_{H'}(v) + |\alpha^{-1}(i) \cap N^-(v)| <   {f_i'}^-(v) + |\alpha^{-1}(i) \cap N^-(v)| = f_i^-(v). \]
    Symmetrically, $d^+_{H'}(v)<{f_i'}^+(v)$ implies $d^+_{H}(v)<f^+_i(v)$. We are thus done.   
    
    \medskip
    
    The algorithm is elementary. It consists on a loop over vertices $v\in X$, updating the table storing the colouring with colour $\alpha(v)$ for $v$, deleting $v$ from $D$, and visiting every neighbour $u$ of $v$ in order to update its adjacency list (that is removing $v$ from $N^-(u)$ and/or $N^+(u)$), and its colour budget $f_{\alpha(v)}(u)$. This is clearly linear in $\sum_{v\in X}d(v)$.
\end{proof}

Let $(D,F)$ be a valid pair.
Let $X$ be a subset of vertices of $D$, $\alpha : X \to [s]$ be a partial colouring of $D$, and $(D',F')$ be the pair reduced from $(D,F)$ by $\alpha$.
We say that the colouring $\alpha$ of $D\ind{X}$ is \textit{safe} if each of the following holds:
\begin{itemize}
    \item $\alpha$ is an $F$-dicolouring of $D\ind{X}$,
    \item $D'$ is connected, and
    \item $(D,F)$ is a hard pair if and only if $(D',F')$ is a hard pair. 
\end{itemize}
For the particular case $X=\{v\}$, we thus say that colouring $v$ with $c$ is \textit{safe} if $f_c(v) \neq (0,0)$, $D-x$ is connected, and the reduced pair $(D',F')$ is a hard pair if and only if $(D,F)$ is.

\medskip

The following lemma tells us how variable bidegeneracy relates to vertex orderings.

\begin{lemma}\label{lemma:equiv-var-bidegeneracy}
    Given a digraph $D=(V,A)$, and a function $f : V \xrightarrow{} \NN^2$,
    $D$ is strictly-$f$-bidegenerate
    if and only if there exists an ordering $v_1,\ldots,v_n$ of the vertices of $D$ such that
    \begin{equation*}
        f^-(v_i) > |N^-(v_i)\cap \{v_j\ |\ j\leq i\}| \mbox{~~~~or~~~~}  f^+(v_i) > |N^+(v_i)\cap \{v_j\ |\ j\leq i\}|
    \end{equation*}
\end{lemma}
\begin{proof}
    $(\Longrightarrow)$ We proceed by induction on the number of vertices $n$. This implication clearly holds for $n=1$. Assume now that $n\geq 2$. Let $v_n$ be a vertex of $D$ such that $f^-(v_n) > d^-(v_n)$ or $f^+(v_n) > d^+(v_n)$, which exists by strict bidegeneracy of $D$.
    Since $d^-(v_n)\geq |N^-(v_n)\cap S|$ and $d^+(v_n)\geq |N^+(v_n)\cap S|$ for any set $S \subseteq V$, the property holds for $v_n$.
    By induction there is an ordering $v_1,\ldots,v_{n-1}$ of the vertices of $D'=D - v_n$ such that for every vertex $v_i$, with $1\leq i\leq n-1$, we have
    \begin{equation*}
        f^-(v_i) > |N_{D'}^-(v_i)\cap \{v_j\ |\ j\leq i\}| \mbox{~~~~or~~~~}  f^+(v_i) > |N_{D'}^+(v_i)\cap \{v_j\ |\ j\leq i\}|
    \end{equation*}
    As $N_{D'}^-(v)\cap \{v_j\ |\ j\leq i\} = N^-_D(v)\cap \{v_j\ |\ j\leq i\}$ and $N_{D'}^+(v)\cap \{v_j\ |\ j\leq i\} = N^+_D(v)\cap \{v_j\ |\ j\leq i\}$, the property still holds in $D$ and we are done.
    
    $(\Longleftarrow)$ For any subdigraph $H=(V_H,A_H)$ of $D$, let $v_k$ be the vertex of $V_H$ with largest index $k$. Since $N^-_H(v_k) = N^-_D(v_k) \cap V_H \subseteq N^-_D(v_k)\cap \{v_j\ |\ j\leq k\}$, and $N^+_H(v_k) = N^+_D(v_k) \cap V_H \subseteq N^+_D(v_k)\cap \{v_j\ |\ j\leq k\}$, we have that 
    \[
    f^-(v_k) > |N^-_D(v)\cap \{v_j\ |\ j\leq k\}|\geq d_H^-(v_k) \mbox{~~~~or~~~~}  f^+(v_i) > |N^+_D(v)\cap \{v_j\ |\ j\leq i\}|\geq d_H^+(v_k).
    \]
\end{proof}

\subsection{Greedy algorithm}\label{subsec:greedy-algorithm}

In this subsection we consider an algorithm that greedily colours $D$, partially or entirely.
We are going to give conditions ensuring that this approach succeeds in providing a (partial) $F$-dicolouring. 

\begin{lemma}\label{lemma:validity-algo-greedy}
    Given a (non-necessarily valid) pair $(D,F)$ and an ordered list of vertices $v_1,\ldots,v_\ell$ of $V(D)$, there is an  algorithm, running in time $O\left(\sum_{i=1}^\ell d(v_i)\right)$, that tries to colour the vertices $v_1,\ldots,v_\ell$ in order to get an $F$-dicolouring
    of $D[v_1,...,v_\ell]$.
    If the algorithm succeeds, it also computes the reduced pair $(D',F')$.
    
    If $F$ is such that for every $v_i\in\{v_1,\ldots,v_\ell\}$, we have
    \begin{equation}
        \sum_{c=1}^s f_c^-(v_i) > |N^-(v_i) \cap \{v_j \ |\ j \leq i\}| \mbox{~~~~or~~~~} \sum_{c=1}^s f_c^+(v_i) > |N^+(v_i)\cap \{v_j \ |\ j \leq i\}|\tag{$\star \star$}
        \label{property_degenerate_F}
    \end{equation}
    the algorithm succeeds. If condition~\eqref{property_degenerate_F} is fulfilled, if $(D,F)$ is valid, and if $D-\{v_1,\ldots,v_\ell\}$ is connected, the reduced pair $(D',F')$ is a valid pair.
\end{lemma}

Note that as $\sum_{v\in D} d(v)=2|A(D)|$ the time complexity here is $O(|V(D)|+|A(D)|)$, even when the whole digraph is coloured. 

\begin{proof}
    The algorithm simply consists in considering the vertices $v_1,\ldots,v_\ell$ in this order and, if possible, to colour $v_i$ with a colour $c$ such that  $f^-_c(v_i)> \left|\left\{ u\in N^-(v_i)\cap \{v_1,\ldots,v_i\}\ |\ u \text{ is coloured }c\right\}\right|$, or such that $f^+_c(v_i)> \left|\left\{ u\in N^+(v_i)\cap \{v_1,\ldots,v_i\}\ |\ u \text{ is coloured }c\right\}\right|$.
    
    The complexity of the algorithm holds because, it actually consists in a loop over vertices $v_1,\ldots,v_\ell$, where 1) it looks for a colour $c$ such that $f_c(v_i)\neq (0,0)$ (where $F$ is updated after each vertex colouring), and if it finds such a colour, 2) colours $v_i$ and updates $(D,F)$ into the corresponding reduced pair $(D',F')$. Step 1) is done in constant time (see Subsection~\ref{subsec:data_structure}), and step 2) is done in $O(d(v_i))$ time (by Lemma~\ref{lemma:safe-colouring}).
    Hence, the complexity clearly follows. Lemma~\ref{lemma:equiv-var-bidegeneracy}, applied to $D\ind{ \{v_1,...,v_\ell\} }$, implies that if the algorithm succeeds, the obtained colouring is an $F$-dicolouring of $D\ind{ \{v_1,...,v_\ell\} }$.
    
    To show the second statement, let us show the following invariant of the algorithm: 
    
    \begin{center}
        \begin{minipage}{0.9\textwidth}
            {\it 
                $(\clubsuit)$ At the beginning of the $i^\text{th}$ iteration of the main loop, for any vertex $v_k$ with $i\leq k\leq \ell$, at least one of the following occurs:
                \begin{itemize}[leftmargin=*]
                    \item its number of uncoloured in-neighbours with index at most $k$ is less than $\sum_{c=1}^s f_c^-(v_i)$, or
                    \item its number of uncoloured out-neighbours with index at most $k$ is less than $\sum_{c=1}^s f_c^+(v_i)$,
                \end{itemize}
                where the sums are made on the updated functions $f_c$. 
            }
        \end{minipage}
    \end{center}  
    
    Note that, by~\eqref{property_degenerate_F}, $(\clubsuit)$ holds for the first iteration.
    Let us show that if $(\clubsuit)$ holds at the beginning of the $i^\text{th}$ iteration for some vertex $v_k$ with $i<k$, then it still holds at the beginning of the $i+1^\text{th}$ iteration.
    Indeed, during the $i^\text{th}$ iteration, if the sum $\sum_{c=1}^s f_c^-(v_k)$ decreases, it decreases by exactly one, and in that case we have $v_i\in N^-(v_k)$. Hence, the number of uncoloured in-neighbours of $v_k$ with index at most $k$ also decreases by one, and $(\clubsuit)$ sill holds. The same holds for the out-neighbourhood of $v_i$.
    
    By $(\clubsuit)$, at the beginning of the $i^\text{th}$ iteration, we have $\sum_{i=1}^s f_i(v) \neq (0,0)$. Hence, there is a colour $c$ such that $f_c(v_i)\neq (0,0)$, so there is always a colour (\textit{e.g.} $c$) available for colouring $v_i$, and the algorithm thus succeeds in colouring the whole digraph $D[v_1,\ldots,v_\ell]$. 
    
    Finally for the last statement, if~\eqref{property_degenerate_F} holds, then the algorithm succeeds in producing an $F$-dicolouring.
    Furthermore, if $(D,F)$ is valid, for any vertex $u\in V(D)\setminus \{v_1,\ldots,v_\ell\}$, we have $d_D^-(u)\leq \sum_{c=1}^s f_c^-(v_i)$ and $d_D^+(u)\leq \sum_{c=1}^s f_c^+(v_i)$. For each of these inequalities, and at each iteration of the main loop, the left hand side decreases if and only if the right hand side does, in the reduced pair $(D',F')$. Hence,~\eqref{property_D_F} holds in $(D',F')$, and since $D'$ is connected, the pair $(D',F')$ is valid.
\end{proof}

\subsection{Solving loose instances}\label{subsec:reduction-tightening}

If the input digraph $D$ as a whole is strictly-$\Tilde{f}$-bidegenerate for some $\Tilde{f}$, it may be easy to produce an $F$-dicolouring, under some conditions on $\Tilde{f}$.
Note that, in what follows, we do not ask for $(D,F)$ to be a valid pair, so $D$ may be disconnected and vertices $v$ do not necessarily satisfy~\eqref{property_D_F}.

\begin{lemma}
    \label{lemma:partition_degenerate}
    Let $D=(V,A)$ be a digraph, $F=(f_1,\ldots,f_s)$ be a sequence of functions $f_i : V \to \mathbb{N}^2$, and $\Tilde{f}=\sum_{i=1}^s f_i$. If $D$ is strictly-$\Tilde{f}$-bidegenerate, then $D$ is $F$-dicolourable and an $F$-dicolouring can be computed in linear time.
\end{lemma}
\begin{proof}
    Let $(D=(V,A),F=(f_1,\ldots,f_s))$ be such a pair. By Lemma~\ref{lemma:equiv-var-bidegeneracy}, there exists an ordering $\sigma=v_1,\ldots,v_n$ of $V(D)$ such that, for every $i\in[s]$, $v_i$ satisfies $d^+_{D_i}(v_i) < f^+(v_i)$ or $d^-_{D_i}(v_i) < f^-(v_i)$, where $D_i$ is the subdigraph of $D$ induced by $\{v_1,\ldots,v_i\}$.
    
    \medskip
    
    We prove that such an ordering $\sigma$ can be computed in linear time. Our proof  follows a classical linear-time algorithm for computing the usual degeneracy of a graph, but we give the proof for completeness. We create two tables $out\_gap$ and $in\_gap$, both of size $n$ and indexed by $V$. We iterate once over the vertices of $V$ and, for every vertex $v$, we set $out\_gap(v)$ to $\Tilde{f}^+(v) - d^+(v)$ and $in\_gap$ to $\Tilde{f}^-(v) - d^-(v)$. We also store in a set $S$ all vertices $v$ for which $out\_gap(v)\leq -1$ or $in\_gap(v)\leq -1$.
    
    Then for $i$ going from $n$ to $1$, we choose a vertex $u$ in $S$ that has not been treated before, we set $v_i$ to $u$, and for every in-neighbour (respectively out-neighbour) $w$ of $u$ that has not been treated before, we decrease $in\_gap(w)$ (respectively $out\_gap(w)$) by one. If $out\_gap(w)\leq -1$ or $in\_gap(w) \leq -1$, we add $w$ to $S$.
    We then remember (using a boolean table for instance) that $u$ has been treated. 
    
    Following this linear-time algorithm, at the beginning of each step $i\in [n]$, for every non-treated vertex $u$, we have $in\_gap = f^-(u) - d^-_{D_i}(u)$ and $out\_gap = f^+(u) - d^+_{D_i}(u)$. We also maintain that $S$ contains all the vertices $u$ satisfying  $in\_gap(u) \leq -1$ or $out\_gap(u)\leq -1$.
    Since $D$ is strictly-$\Tilde{f}$-bidegenerate, at step $i$, there exists a non-treated vertex $u\in S$.
    
    \medskip
    
    We now admit that such an ordering $\sigma=(v_1,\ldots,v_n)$ has been computed, and we greedily colour the vertices from $v_1$ to $v_n$. The result follows from Lemma~\ref{lemma:validity-algo-greedy} and the fact that $\Tilde{f} = \sum_{i=1}^s f_i$.
\end{proof}

A valid pair is \textit{tight}, if for every vertex the two inequalities of~\eqref{property_D_F} are equalities.
The following particular case of Lemma~\ref{lemma:partition_degenerate} treats the case of instances that are non-tight, which we also refer to as \textit{loose}.

\begin{lemma}
    \label{lemma:partition_loose}
    Let $(D,F)$ be a valid pair that is loose. Then $D$ is $F$-dicolourable and an $F$-dicolouring can be computed in linear time.
\end{lemma}
\begin{proof}
    Let $(D=(V,A),F=(f_1,\ldots,f_s))$ be such a pair and let $\Tilde{f} : V \to \mathbb{N}^2$ be $\sum_{i=1}^s f_i$. It is sufficient to prove that $D$ is necessarily strictly-$\Tilde{f}$-bidegenerate, so the result follows from Lemma~\ref{lemma:partition_degenerate}.
    
    Assume for a contradiction that $D$ is not strictly-$\Tilde{f}$-bidegenerate, so by definition there exists an induced subdigraph $H$ of $D$ such that, for every vertex $v\in V(H)$, $d^+_H(v) \geq \Tilde{f}^+(v)$ and $d^-_H(v) \geq \Tilde{f}^-(v)$. By definition of $\Tilde{f}$ and because $(D,F)$ is a valid pair, we thus have $d^+_H(v) = \sum_{i=1}^sf_i^+(v)$ and $d^-_H(v) = \Tilde{f}^-(v)$. We directly deduce that $H\neq D$ as $(D,F)$ is loose.
    
    Since $D$ is connected (as $(D,F)$ is a valid pair) and because $H\neq D$, there exists in $D$ an arc between vertices $u$ and $v$ such that $u\in V(H)$ and $v\in V(D)\setminus V(H)$. Depending on the orientation of this arc, we have $d^+_H(u) < d^+_D(u) = \sum_{i=1}^sf_i^+(u)$ or $d^-_H(u) < d^-_D(u) = \sum_{i=1}^sf_i^-(u)$, a contradiction.
\end{proof}

\subsection{Reducing to a block and detecting hard pairs}\label{subsec:reduction-block}

We have just shown how to solve loose instances.
In this subsection, we thus consider a tight instance $(D,F)$, and show how to obtain a reduced instance $(B,F')$ where $B$ is a block of $D$, by safely colouring $V(D) \setminus V(B)$. Then, $(B,F')$ is a hard pair if and only if $(D,F)$ is, in which case we may terminate the algorithm. Otherwise, the following subsections show that $B$ may be $F'$-dicoloured, and together with the colouring of $V(D) \setminus V(B)$ fixed at this step, this yields an $F$-dicolouring of $D$.

Our reduction of $(D,F)$ proceeds by considering the end-blocks of $D$ one after the other. If an end-block $B$, together with $F$, may correspond to a monochromatic hard pair, a bicycle hard pair, or a complete hard pair glued to the rest of the digraph, then we safely colour $V(B) \setminus V(D)$ and move to the next end-block. If $B$ cannot be such a block, then we will show that we can safely colour $V(D) \setminus V(B)$.

Before formalising this strategy, we need a few definitions. Given a valid pair $(D,F)$, an end-block $B$ of $D$, with cut-vertex $x$, is a \textit{hard end-block} if it is of one of the following types:
\begin{enumerate}[label=$(\roman*)$]
    \item There exists a colour $i\in [s]$ such that $f_i(x) \geq (d^-_B(x),d^+_B(x))$, and for every $v\in V(B)\setminus\{x\}$ we have $f_i(v) = (d^-_B(v),d^+_B(v))$ and $f_k(v) = (0,0)$ when $k\neq i$.
    
    We refer to such a hard end-block as a \textit{monochromatic hard end-block}.
    
    \item $B$ is a bidirected odd cycle and there exists colours $i\neq j$ such that, $f_i(x) \geq (1,1)$, $f_j(x)\geq (1,1)$ and for every $v\in V(B)\setminus\{x\}$, we have $f_c(v) = (1,1)$ if $c\in \{i,j\}$ and $f_c(v) = (0,0)$ otherwise. 
    
    We refer to such a hard end-block as a \textit{bicycle hard end-block}.
    
    \item $B$ is a bidirected complete graph, the functions $f_1,\ldots,f_s$ are constant and symmetric on $V(B)\setminus \{x\}$, and $f_i(x) \geq f_i(u)$ for every $i\in [s]$ and every $u\in V(B) \setminus \{x\}$.
    
    We refer to such a hard end-block as a \textit{complete hard end-block}.
\end{enumerate}

Let $(D,F=(f_1,\ldots,f_s))$ be a tight valid pair such that $D$ is not biconnected and let $B$ be a hard end-block of $D$ with cut-vertex $x$. Let $u$ be any vertex of $B$ distinct from $x$. We define the \textit{contraction} $(D',F'=(f_1',\ldots,f_s'))$ of $(D,F)$ with respect to $B$, as follows: 
\begin{itemize}
    \item $D' = D- (V(B) \setminus \{x\})$;
    \item for every vertex $v\in V(D')\setminus \{x\}$ and every $i\in [s]$, $f_i'(v) = f_i(v)$;
    \item if $B$ is a monochromatic hard end-block, let $c$ be the unique colour such that $f_c(u) \neq (0,0)$, then $f_c'(x) = f_c(x) - (d^-_B(x),d^+_B(x))$ and $f_i'(x) = f_i(x)$ for every $i\in [s]\setminus \{c\}$;
    \item otherwise, $B$ a bicycle or a complete hard end-block, and $f_i'(x) = f_i(x) - f_i(u)$ for every $i\in [s]$.
\end{itemize}

\begin{lemma}~\label{lemma:equiv-hard-pair}
    Let $(D,F=(f_1,\ldots,f_s))$ be a tight valid pair such that $D$ is not biconnected and let $B$ be a hard end-block of $D$ with cut-vertex $x$. Let $(D',F'=(f_1',\ldots,f_s'))$ be the contraction of $(D,F)$ with respect to $B$. Then $(D,F)$ is a hard pair if and only if $(D',F')$ is a hard pair.
\end{lemma}

\begin{proof}
    If $(D',F')$ is a hard pair, $(D,F)$ could be defined as the join hard pair obtained from two hard pairs, $(D',F')$ and $(B,\Tilde{F}=(\Tilde{f}_1,\ldots,\Tilde{f}_s))$, where $\Tilde{f}_i = f_i - f_i'$ for every $i\in [s]$.
    
    Conversely, let us show by induction that for every hard pair $(D,F)$, and every end-block $B$, then $B$ is a hard end-block and $(D',F')$, the contraction of $(D,F)$ with respect to $B$, is a hard pair.
    This is trivial if $D$ is biconnected as $D$ does not have any end-block. We thus assume that $(D,F)$ is a join hard pair obtained from two hard pairs $(D_1,F_1)$ and $(D_2,F_2)$. Assume without loss of generality that $B$ is a block of $D_1$. If $B$ is not an end-block in $D_1$, since it is an end-block of $D$, we necessarily have $D_1 = B$. Hence $(D',F')$ is exactly $(D_2,F_2)$, which is a hard pair. Furthermore, irrespective of the type of hard pair $(D_1,F_1)$, by the definition of hard join, $D_1 = B$ is clearly a hard end-block of $D$.
    If $B$ is an end-block of $D_1$, By induction hypothesis, $B$ is a hard end-block of $(D_1,F_1)$ and the 
    contraction of $(D_1,F_1)$ with respect to $B$, $(D'_1,F'_1)$, is a hard pair. It is easy to check that $(D',F')$ is exactly the join hard pair obtained from $(D'_1,F'_1)$ and $(D_2,F_2)$. 
\end{proof}

Before describing our algorithm reducing the instance to a single block, we need the following subroutine testing if an end-block is hard.

\begin{lemma}
    \label{lemma:check_hard-endblock}
    Given a tight valid pair $(D,F=(f_1,\ldots,f_s))$, and an end-block $B$ of $D$ with cut-vertex $x$, testing whether $B$ is a hard end-block can be done in time $O(|V(B)|+|A(B)|)$.
\end{lemma}

\begin{proof}
    The algorithm takes block $B$ with its cut-vertex $x$, and considers any $u \in V(B) \setminus \{x\}$ as a reference vertex. 
    
    We first check whether $B$ is a monochromatic hard end-block. To do so, we first check whether exactly one colour $i$ is available for $u$ (\textit{i.e.} $\{j\mid f_j(u) \neq (0,0)\} = \{i\}$). We then check whether $i$ is indeed the only available colour for every vertex $v\in V(B) \setminus \{u,x\}$. We finally check whether $f_i(x) \geq (d^-_D(x), d^+_D(x))$. Since $(D,F)$ is tight, $B$ is a monochromatic hard end-block if and only if all these conditions are met.
    
    Assume now that $B$ is not a monochromatic hard end-block. 
    Now, $B$ is a hard end-block if and only if it is a bicycle hard end-block or a complete hard end-block. In both cases, the functions $f_1,\ldots,f_s$ must be symmetric, and constant on $V(B) \setminus \{x\}$.    
    Since $(D,F)$ is tight, for every vertex $v \in V(B) \setminus \{x\}$, we can iterate over $\{j \mid f_j(v) \neq (0,0)\}$ in time $O(d(v))$. This allows us to check in linear time whether the functions $f_1,\ldots,f_s$ are symmetric and constant on $V(B) \setminus \{x\}$. If this is not the case,  $B$ is not a hard end-block. We can assume now that $f_1,\ldots,f_s$ are symmetric and constant on $V(B)\setminus \{x\}$. We also assume that $B$ is either a bidirected odd cycle or a bidirected complete graph, for otherwise $B$ is clearly not a hard end-block.
    
    If $B$ is a bidirected odd cycle, we check in constant time that $f_i(u) \neq (0,0)$ holds for exactly two colours. If this does not hold, then $B$ is not a hard end-block. 
    If this is the case, denote these two colours $i,j$, and note (by tightness and symmetry) that $f_i(u)=f_j(u)=(1,1)$, and that $f_k(u)=(0,0)$ for $k\notin \{i,j'\}$. Now $B$ is a bicycle hard end-block if and only if $f_k(x) \geq (1,1)$ for $k\in \{i,j\}$, which can be checked in constant time.
    
    Assume finally that $B$ is a bidirected complete graph. Since $(D,F)$ is tight, and because the functions $f_1,\ldots,f_s$ are symmetric and constant on $V(B) \setminus \{x\}$, then $B$ is a complete hard end-block if and only if $f_i(x) \geq f_i(u)$ for every colour $i \in [s]$. This can be done in $O(|V(B)|)$ time, since there are at most $|V(B)|$ colours $i$ such that $f_i(u)\neq (0,0)$ (as $(D,F)$ is tight).
\end{proof}

We now turn to showing the main lemma of this subsection, which reduces any tight valid pair to a biconnected one.
\begin{lemma}
    \label{lemma:reduction_to_block}
    Let $(D,F)$ be a tight valid pair. There exists a block $B$ of $D$ such that $V(D) \setminus V(B)$ may be safely coloured, yielding the reduced pair $(B,F')$. Moreover, $(B,F')$ and the colouring of $V(D) \setminus V(B)$ can be computed in linear time.
\end{lemma}
\begin{proof}
    Let $(D,F)$ be such a pair. We first compute, in linear time~\cite{tarjanSIC14}, an ordering $B^1,\ldots,B^r$ of the blocks of $D$ and an ordered set of vertices $(x_2,\ldots,x_r)$ in such a way that for every $2\leq \ell \leq r$, $x_\ell$ is the only cut-vertex of $B^\ell$ in $D^\ell = D\ind{ \bigcup_{j=1}^\ell V(B^\ell)}$. If $D$ is biconnected, that is $r=1$, the result follows for $B=D$ and there is nothing to do. We now assume $r\geq 2$.
    
    For $\ell$ going from $r$ to $2$, we proceed as follows. We consider the block $B^\ell$, which is an end-block of $D^\ell$. We will either safely colour $V(D^\ell) \setminus V(B^\ell)$ and output the reduced pair $(B^\ell,F')$, or safely colour $V(B^\ell)\setminus\{x_\ell\}$, compute the reduced pair $(D^{\ell-1},F^{\ell-1})$, and go on with the end-block $B^{\ell-1}$ of $D^{\ell-1}$. If we colour all the blocks $B^\ell$ with $\ell \geq 2$, we output the reduced pair $(D^1,F^1)$.
    
    More precisely, when considering $B^\ell$, we first check in time $O(|V(B^\ell)|+|A(B^\ell)|)$ whether $B^\ell$ is a hard end-block of $(D^\ell,F^\ell=(f_1^\ell,\ldots,f_s^\ell))$ (by Lemma~\ref{lemma:check_hard-endblock}).
    We distinguish two cases, depending on whether $B^\ell$ is a hard end-block.
    \begin{description}
        \item[Case 1:] \textit{$B^\ell$ is a hard end-block of $D^\ell$}.
        
        In this case, we safely colour the vertices of $V(B^\ell)\setminus \{x_\ell\}$ in time $O(|V(B^\ell)|+|A(B^\ell)|)$.
        To do so, we first compute an ordering $\sigma = (v_1,\ldots,v_b=x_\ell)$ of $V(B^\ell)$ corresponding to a leaves-to-root ordering of a spanning tree of $B^\ell$ rooted in $x_\ell$, in time $O(|V(B^\ell)|+|A(B^\ell)|)$. 
        We then greedily colour the vertices $v_1,\ldots,v_{b-1}$ in this order.
        For every $j\in [b-1]$, the vertex $v_j$ has at least one neighbour in $\{v_{j+1},\ldots,v_b\}$ (its parent in the spanning tree), so condition~\eqref{property_degenerate_F} is fulfilled. Hence, Lemma~\ref{lemma:validity-algo-greedy} ensures that the algorithm succeeds and that the reduced pair $(D^{\ell-1},F^{\ell-1}=(f_1^{\ell-1},\ldots,f_s^{\ell-1}))$ is a valid pair. 
        We will now show that $(D^{\ell-1},F^{\ell-1})$ is necessarily tight, and that the colouring of $V(B^\ell)\setminus \{x_\ell\}$ we computed is safe, by showing that $(D^{\ell-1},F^{\ell-1})$ is exactly the contraction of $(D^{\ell},F^{\ell})$ with respect to $B^\ell$.
        
        Assume first that $B^\ell$ is a monochromatic hard end-block. Then, all the vertices of $V(B)\setminus \{x_\ell\}$ are coloured with the same colour $i$, and by definition of a monochromatic hard end-block, we have $f^{\ell-1}_i(x_\ell) = f^{\ell}_i(x_\ell) - (d^-_{B^\ell}(x_\ell), d^+_{B^\ell}(x_\ell))$ and $f^{\ell-1}_{j}(x_\ell) = f^{\ell}_{j}(x_\ell)$ for every colour $j\neq i$. Hence, $(D^{\ell-1},F^{\ell-1})$ is tight (as the degrees from $D^\ell$ to $D^{\ell-1}$ only change for $x_\ell$, for which they decrease exactly by $d^-_{B^\ell}(x_\ell)$, and $d^+_{B^\ell}(x_\ell)$) and Lemma~\ref{lemma:equiv-hard-pair} implies that the colouring we computed is safe.
        
        Assume now that $B^\ell$ is a bicycle hard end-block. Let $i \neq j$ be such that $f_{i}^\ell(u) = f_{j}^\ell(u) = (1,1)$ for every vertex $u\in V(B^\ell) \setminus \{x_\ell\}$. Recall that, by definition of a bicycle hard end-block, $f_{k}^\ell (x_\ell) \geq (1,1)$ for $k\in \{i,j\}$. Let $u,v$ be the two neighbours of $x_\ell$ in $B^\ell$, then $u$ and $v$ are connected by a bidirected path of odd length in $B^\ell - x_\ell$. This implies that $u$ and $v$ are coloured differently, and $f^{\ell-1}_{k}(x_\ell) = f^{\ell}_k(x_\ell) - (1,1)$ for $k\in \{i,j\}$ and $f^{\ell-1}_{k}(x_\ell) = f^{\ell}_{k}(x_\ell)$ for $k\notin \{i,j\}$. Since $d^+_{B^\ell}(x_\ell) = d^-_{B^\ell}(x_\ell) = 2$, we obtain that $(D^{\ell-1},F^{\ell-1})$ is tight, and Lemma~\ref{lemma:equiv-hard-pair} implies that the colouring we computed is safe.
        
        Assume finally that $B^\ell$ is a complete hard end-block, and let $u$ be any neighbour of $x_\ell$ in $B$. The functions $f_1^\ell,\ldots,f_s^\ell$ are constant and symmetric on $V(B) \setminus \{x_\ell\}$. We also have $f_i^\ell(x_\ell)\geq f_i^\ell(u)$ for every $i\in [s]$.
        By construction of the greedy colouring, for every $i\in [s]$, at most ${f_i^\ell}^+(u)$ vertices of $V(B^\ell) \setminus \{x_\ell\}$ are coloured $i$. Since $|V(B^\ell)|-1$ vertices are coloured in total, and because $\sum_{i=1}^s {f_i^\ell}^+(u) = |V(B^\ell)| -1$, we conclude that, for every $i\in [s]$, exactly ${f_i^\ell}^+(u)$ vertices of $V(B^\ell) \setminus \{x_\ell\}$ are coloured $i$. Hence, we obtain that $(D^{\ell-1},F^{\ell-1})$ is tight, and Lemma~\ref{lemma:equiv-hard-pair} implies that the colouring we computed is safe.
        
        \item[Case 2:] \textit{$B^\ell$} is not a hard end-block of $D^\ell$.
        
        In this case, we colour the vertices of $V(D^{\ell-1})\setminus \{x_\ell\}$ in time $O(|V(D)|+|A(D)|)$ as follows. We first compute an ordering $\sigma = (v_1,\ldots,v_{n'}=x_\ell)$ of $V(D^{\ell-1})$, that is a leaves-to-root ordering of a spanning tree of $D^{\ell-1}$ rooted in $x_\ell$, in time $O(|V(D^\ell)|+|A(D^\ell)|)$. We greedily colour vertices $v_1,\ldots,v_{n'-1}$ in this order. For every $j\in [n'-1]$, the vertex $v_j$ has at least one neighbour in $\{v_{j+1},\ldots,v_{n'}\}$ (its parent in the spanning tree), so condition~\eqref{property_degenerate_F} is fulfilled. Hence by Lemma~\ref{lemma:validity-algo-greedy}, the greedy colouring succeeds and provides a reduced pair $(B^{\ell},F'=(f_1',\ldots,f_s'))$ that is valid. 
        
        Finally observe that $(B^\ell,F')$ is not a hard pair since, for every $i\in [s]$, we have $f_i'(x_\ell) \leq f_i^\ell(x_\ell)$. So if $(B^\ell,F')$ is a hard pair, $B^\ell$ is necessarily a hard end-block of $(D^\ell,F^\ell)$, a contradiction. The result follows.
    \end{description}
    
    Note that Case~2 may only be reached once, at which point it outputs a reduced pair.
    Therefore, the running time of the algorithm described above is bounded by $O(\sum_\ell \left(|V(B^\ell)|+|A(B^\ell)|\right) + |V(D)|+|A(D)|)$, which is linear in $|V(D)|+|A(D)|$.
    Assume finally that, in the process above, the second case is never attained. Then the result follows as $B^1$ is a block of $D$, and we found a safe colouring of $V(D) \setminus V(B^1)$.
\end{proof}

With the last two lemmas at our disposal, we are ready to test whether a tight valid pair is hard.
\begin{lemma}
    \label{lemma:check_hard-pair}
    Given a tight valid pair $(D,F=(f_1,\ldots,f_s))$, testing whether $(D,F)$ is a hard pair can be done in linear time.
\end{lemma}
\begin{proof}
    We first consider the pair $(B,F')$ reduced from $(D,F)$, where $B$ is a block of $D$, obtained by safely colouring $V(D) \setminus V(B)$ through Lemma~\ref{lemma:reduction_to_block}.
    In particular, $(B,F')$ is a hard pair if and only if $(D,F)$ is, and as $B$ is biconnected, we may only check whether $(B,F')$ is a biconnected hard pair, either monochromatic, bicycle or complete.
    This amounts to verifying the conditions of the definition. The fact that this can be done in linear time is already justified in Lemma~\ref{lemma:check_hard-endblock}. Indeed, testing $(B,F')$ for a hard biconnected pair amounts to testing the conditions for a hard end-block, except for the condition on the cut-vertex. Since there is no cut-vertex in $B$, we test the same conditions on all vertices of $B$.
\end{proof}

\subsection{Reducing to pairs with two colours}\label{subsec:reduction-2colblocks}

We have just shown how to reduce any tight valid pair into a biconnected one, after which we may decide whether the initial pair was a hard one.
In this subsection, we therefore consider any tight valid pair $(D,F)$ that is not hard, and such that $D$ is biconnected. We show how the problem of $F$-dicolouring $D$ boils down to $\Tilde{F}$-dicolouring $D$ where $\Tilde{F}$ involves only $s=2$ colours.

\begin{lemma}
    \label{lemma:reduce_to_2colours}
    Let $(D,F=(f_1,\ldots,f_s))$ be a valid pair that is not a hard pair and such that $D$ is biconnected. In linear time, we can either find an $F$-dicolouring of $D$ or compute a valid pair $(D,\Tilde{F}=(\Tilde{f}_1,\Tilde{f}_2))$ such that:
    \begin{itemize}[noitemsep, topsep=0pt]
        \item $(D,\Tilde{F})$ is not hard, and
        \item given an $\Tilde{F}$-dicolouring of $D$, we can compute an $F$-dicolouring of $D$ in linear time.
    \end{itemize}
\end{lemma}
\begin{proof}
    We assume that $(D,F)$ is tight, otherwise we compute an $F$-dicolouring of $D$ in linear time by Lemma~\ref{lemma:partition_loose}.
    Since $(D,F)$ is tight, we can iterate over $\{i \mid f_i(v) \neq (0,0)\}$ in time $O(d(u))$. This allows us to check, in linear time, which of the cases below apply to $(D,F)$.
    
    \begin{description}
        \item[Case 1:] \textit{There exist $u \in V(D)$ and $i \in [s]$ such that $f_i^+(u) \neq f_i^-(u)$.}
        
        We define $\Tilde{f}_1 = f_i$ and $\Tilde{f}_2 = \sum_{j\in [s], j\neq i}f_j$. 
        Let us show that $(D,\Tilde{F})$ is not a hard pair. Since $D$ is biconnected, $(D,\Tilde{F})$ is not a join hard pair. Since $\Tilde{f}_1$ is not symmetric by choice of $\{u,c\}$, $(D,\Tilde{F})$ is neither a bicycle hard pair nor a complete hard pair. Finally, assume for a contradiction that $(D,\Tilde{F})$ is a monochromatic hard pair. Since $\Tilde{f}_1(u)\neq (0,0)$, we thus have $\Tilde{f}_2(v) = (0,0)$ and $f_i(v) = (d^-(v),d^+(v))$ for every vertex $v$. We conclude that $(D,F)$ is also a monochromatic hard pair, a contradiction.
        
        \item[Case 2:] \textit{There exist $u,v\in V(D)$ and $i\in [s]$ such that $f_i(u) \neq f_i(v)$.} 
        
        As in Case~1,  we set $\Tilde{f}_1 = f_i$ and $\Tilde{f}_2 = \sum_{j\in [s], j\neq i}f_j$. Since $D$ is biconnected, $(D,\Tilde{F}=(\Tilde{f}_1,\Tilde{f}_2))$ is not a join hard pair. Moreover, since $f_i(v) \neq f_i(u)$, $\Tilde{f}_1$ is not constant so $(D,\Tilde{F})$ is neither a bicycle hard pair nor a complete hard pair. 
        Again, assume for a contradiction that $(D,\Tilde{F})$ is a monochromatic hard pair. Since $\Tilde{f}_1(u)$ or $\Tilde{f}_1(v)$ is distinct from $(0,0)$, we thus have $\Tilde{f}_2(v) = (0,0)$ and $f_i(v) = (d^-(v),d^+(v))$ for every vertex $v$. We conclude that $(D,F)$ is also a monochromatic hard pair, a contradiction.
        
        \item[Case 3:] \textit{None of the cases above is matched.}
        
        Therefore, for each $i\in[s]$, $f_i$ is a symmetric constant function. Thus, since $(D,F)$ is not a hard pair, and because $(D,F)$ is tight, $D$ is not a bidirected odd cycle or a bidirected complete graph. Let $i\in [s]$ be such that $f_i$ is not the constant function equal to $(0,0)$. We set $\Tilde{f}_1 = f_i$ and $\Tilde{f}_2 = \sum_{j\in [s], j\neq i}f_j$. Since $(D,F)$ is not a monochromatic hard pair, there is a colour $k\neq i$ such that $f_{k}$ is not the constant function equal to $(0,0)$. Hence, none of $\Tilde{f}_1,\Tilde{f}_2$ is the constant function equal to $(0,0)$, and $(D,(\Tilde{f}_1,\Tilde{f}_2))$ is not a hard pair.
    \end{description}
    
    In each case, we have built a valid pair $(D,\Tilde{F}=(\Tilde{f}_1,\Tilde{f}_2))$ where $\Tilde{f}_1 = f_i$ and $\Tilde{f}_2 = \sum_{j\in [s], j\neq i}f_j$, in such a way that $(D,\Tilde{F})$ is not hard.
    We finally prove that, given an $\Tilde{F}$-dicolouring of $D$, we can compute an $F$-dicolouring of $D$ in linear time. Let $\Tilde{\alpha}$ be an $\Tilde{F}$-dicolouring of $D$. Let $X$ be the set of vertices coloured $1$ in $\Tilde{\alpha}$ and $\Hat{D}$ be $D - X$ ($\hat{D}$ is built in linear time by successively removing vertices in $X$ from $D$). We define $\Hat{F}=(\hat{f}_1,\ldots,\hat{f}_s)$, $\hat{f}_i : V(\Hat{D}) \to \mathbb{N}^2$ as follows:
    \[
    \hat{f}_k(v) = \left\{
    \begin{array}{ll}
        f_k(v) & \mbox{if } k \neq i \\
        (0,0) & \mbox{otherwise.}
    \end{array}
    \right.   
    \]
    
    Since $\Tilde{f}_2$ is exactly $\sum_{i=1}^s\hat{f}_i$, by definition of $\Tilde{\alpha}$ and $\hat{D}$, we know that $\hat{D}$ is strictly-$(\sum_{i=1}^s\hat{f}_i)$-bidegenerate. Then applying Lemma~\ref{lemma:partition_degenerate} yields an $\hat{F}$-dicolouring $\hat{\alpha}$ of $\hat{D}$ in linear time. The colouring $\alpha$ defined as follows is thus an $F$-dicolouring of $D$ obtained in linear time:
    \[
    \alpha(v) = \left\{
    \begin{array}{ll}
        i & \mbox{if } \Tilde{\alpha}(v) = 1 \\
        \hat{\alpha}(v) & \mbox{otherwise. }
    \end{array}
    \right.   
    \]
\end{proof}

\subsection{Solving blocks with two colours - particular cases}\label{subsec:algo-block-2col-particular}

We are now left to deal with pairs $(D,F)$ such that $D$ is biconnected, and $F$ only involves two colours.
Eventually, our strategy consists in finding a suitable decomposition for $D$, and colouring parts of it inductively, which will be done in Subsection~\ref{subsec:algo-block-2col-general}.
In this subsection, we first deal with some specific forms $(D,F)$ may take, showing $D$ can be $F$-dicoloured in those cases.
Along the way, this allows us to deal with increasingly restricted instances, later enabling us to constrain the instances considered in the induction.
In Lemma~\ref{lemma:partition_deficiency} and Lemma~\ref{lemma:partition_unbalanced}, we exhibit an $F$-dicolouring of $D$ if there exists a vertex $x$ which does not satisfy some conditions relating its neighbourhood and $F$.
Then, we solve instances where $D$ is a bidirected complete graph in Lemma~\ref{lemma:partition_complete}, or an orientation of a cycle in Lemma~\ref{lemma:partition_cycle}.
Lastly, we solve those that are constructed by a star attached to a cycle in Lemma~\ref{lemma:partition_wheel}.
These will serve as base cases for the induction.

Given a valid pair $(D,F=(f_1,f_2))$ that is tight, non-hard, and such that $D$ is biconnected, let us define a set of properties that $(D,F)$ may or may not fulfil.
Each of these properties is easy to check in linear time, and we will see in the following that when they are not met, there is a linear-time algorithm providing an $F$-dicolouring of $D$. The first property,~\eqref{prop:deficiency}, guarantees that both $f_1$ and $f_2$ \emph{exceed} some lower bound.
It states that both $f_1^-,f_2^-$ (respectively $f_1^+,f_2^+$) are non-zero on vertices having at least one in-neighbour (respectively out-neighbour).
\begin{equation}
    \forall  x,c~~~\text{ we have }~~~ (d^-(x) > 0 \Longrightarrow f_c^-(x)\geq 1) ~~~\text{ and }~~~(d^+(x) > 0 \Longrightarrow f_c^+(x)\geq 1).
    \tag{E}
    \label{prop:deficiency}
\end{equation}
Given a digraph $D=(V,A)$, we define the function $\mathbb{1}_{A}: V \times V \to \mathbb{N}$ as follows:
\[
\mathbb{1}_{A}(u,v) = \left\{
\begin{array}{ll}
    1 & \mbox{if } uv \in A \\
    0 & \mbox{otherwise.} 
\end{array}
\right.
\]
Note in particular, that if~\eqref{prop:deficiency} holds, then the following holds.
\begin{equation}
    \forall  x,y,z,c~~~\text{ we have }~~~ f_c(x) \geq (\mathbb{1}_{A}(y,x),\mathbb{1}_{A}(x,z)).
    \tag{E'}
    \label{prop:deficiency-prime}
\end{equation}

\begin{lemma}\label{lemma:partition_deficiency}
    Let $(D,F=(f_1,f_2))$ be a valid tight non-hard pair such that $D$ is biconnected. If the pair does not fulfill property~\eqref{prop:deficiency}, then $D$ is $F$-dicolourable. Furthermore, there is a linear-time algorithm checking property~\eqref{prop:deficiency},  and if the property is not met, that computes an $F$-dicolouring of $D$.
\end{lemma}
\begin{proof}
    We first prove that, for some arc $uv$  and for some colour $c\in \{1,2\}$, we have:
    \begin{equation}
        (f_c^+(u)=0 \wedge f_c(v) \neq (0,0)) ~~~\text{ or }~~~(f_c^-(v)=0 \wedge f_c(u) \neq (0,0)).\tag{1}
        \label{ineq:deficient_uv}
    \end{equation}
    
    Note that, if it exists, such an arc is found in linear time. Assume for a contradiction that no such arc exists. By assumption, $\exists x \in V(D),~c\in \{1,2\}$ such that $(d^+(x) > 0 \wedge f_c^+(x) = 0)$ or $(d^-(x) > 0 \wedge f_c^-(x) = 0)$. If $(d^+(x) > 0 \wedge f_c^+(x) = 0)$, let $y$ be any out-neighbour of $x$, then we must have $f_c(y) = (0,0)$ for otherwise $xy$ clearly satisfies~\eqref{ineq:deficient_uv}. 
    Symmetrically, if $(d^-(x) > 0 \wedge f_c^-(x) = 0)$, let $y$ be any in-neighbour of $x$, then we have $f_c(y) = (0,0)$ for otherwise $yx$ satisfies~\eqref{ineq:deficient_uv}. 
    In both cases, we conclude on the existence of a vertex $y$ for which $f_c(y)= (0,0)$. Let $Y$ be the non-empty set of vertices $y$ for which $f_c(y) = (0,0)$ and $Z$ be $V(D) \setminus Y$, that is the set of vertices $z$ for which $f_c(z) \neq (0,0)$.
    Since $(D,F)$ is valid and tight, $Z$ must also be non-empty, for otherwise $(D,F)$ is a monochromatic hard pair. The connectivity of $D$ guarantees the existence of an arc between a vertex in $Y$ and another one in $Z$. This arc satisfies~\eqref{ineq:deficient_uv}.
    
    Once we found an arc $uv$ satisfying~\eqref{ineq:deficient_uv}, we proceed as follows to compute an $F$-dicolouring of $D$. If $f_c^+(u)=0 \wedge f_c(v) \neq (0,0)$, we colour $v$ with $c$, otherwise we have $f_c^-(v)=0 \wedge f_c(u) \neq (0,0)$ and we colour $u$ with $c$. Let $(D',F'=(f_1',f_2'))$ be the pair reduced from this colouring and let $c'$ be the colour distinct from $c$. In the former case, that is $v$ is coloured with $c$, we obtain ${f_{c'}'}^+(u) = f_{c'}^+(u) \geq d^+_D(u) > d^+_{D'}(u)$. In the latter case, $u$ is coloured with $c$ and ${f_{c'}'}^-(v) = f_{c'}^-(v) \geq d^-_D(v) > d^-_{D'}(v)$. In both cases, we obtain that $(D',F')$ is loose, so the result follows from Lemma~\ref{lemma:partition_loose}.
\end{proof}

The next property,~\eqref{prop:digon-sym}, guarantees that the vertices incident only to digons have symmetric constraints.
\begin{equation}
    \forall  x,c~~~\text{ we have }~~~ N^-(x) \neq N^+(x),~~~\text{ or }~~~ f_c^-(x) = f_c^+(x).
    \tag{DS}
    \label{prop:digon-sym}
\end{equation}

Again, we may solve the instance at this point if $(D,F)$ doesn't satisfy~\eqref{prop:digon-sym}.
\begin{lemma}
    \label{lemma:partition_unbalanced}
    Let $(D,F=(f_1,f_2))$ be a valid tight non-hard pair such that $D$ is biconnected. If the pair does not fulfill property~\eqref{prop:digon-sym}, then $D$ is $F$-dicolourable. Furthermore, there is a linear-time algorithm checking property~\eqref{prop:digon-sym},  and if the property is not met, computing an $F$-dicolouring of $D$.
\end{lemma}
\begin{proof}
    Recall that the data structure encoding $D$ allows checking $N^+(u) = N^-(u)$ in constant time (as this is equivalent to checking $|N^+(u)\cap N^-(u)|=d^-(u)=d^+(u)$), so in linear time we may find a vertex $x$ such that $N^-(x) = N^+(x)$ and $f_c^-(x) \neq f_c^+(x)$ for some $c\in \{1,2\}$. 
    
    Let $T$ be a spanning tree rooted in $x$ and let $\sigma=(v_1,\ldots,v_n=x)$ be a leaves-to-root ordering of $V(D)$ with respect to $T$. For every $j\in [n-1]$, the vertex $v_j$ has at least one neighbour in $\{v_{j+1},\ldots,v_{n}\}$ (its parent in the spanning tree), so condition~\eqref{property_degenerate_F} is fulfilled. Hence by Lemma~\ref{lemma:validity-algo-greedy}, there is an algorithm that computes an $F$-dicolouring of $D -x$.
    It remains to show that the reduced pair $(D',F')$ is non-hard. As $D'$ is the single vertex $x$, this is equivalent to showing that $f'_c(x) \neq (0,0)$ for some $c$.
    
    Towards a contradiction, suppose $f'_1(x) = f'_2(x) = (0,0)$.
    For $i\in \{1,2\}$, let $d_i$ be the number of neighbours of $x$ coloured $i$, and note that $d_1+d_2= |N(x)|$. 
    Since $f_i(x) - (d_i,d_i) \leq f'_i(x) = (0,0)$ we have that 
    \begin{equation}
        \forall i~~~~~f_i(x) \leq (d_i,d_i).
        \tag{2}
        \label{eqDS}
    \end{equation} 
    By tightness, we then have that 
    \[(|N(x)|,|N(x)|) = f_1(x) + f_2(x) \leq (d_1+d_2,d_1+d_2) = (|N(x)|,|N(x)|)\]
    Thus, the inequalities of~\eqref{eqDS} are equalities, and both $f_1,f_2$ are symmetric on $x$, a contradiction.
\end{proof}

At this point, we have proven that if one of~\eqref{prop:deficiency} or~\eqref{prop:digon-sym} doesn't hold, an $F$-dicolouring of $D$ may be computed in linear time.
From now on, we thus consider only instances satisfying both conditions, and move on to solving the base cases of our induction, corresponding to $D$ having a certain structure.
We first prove that if $D$ is a bidirected complete graph, and the instance isn't hard, an $F$-dicolouring of $D$ exists and can be computed in linear time.
\begin{lemma}
    \label{lemma:partition_complete}
    Let $(D,F=(f_1,f_2))$ be a valid tight non-hard pair such that $D$ is a bidirected complete graph.
    Then $D$ is $F$-dicolourable and there is a linear-time algorithm providing an $F$-dicolouring of $D$.
\end{lemma}
\begin{proof}
    We assume that $f_1,f_2$ are symmetric ({\it i.e.} $\forall x,c$, $f^+_c(x) = f_c^-(x)$), otherwise we are done by Lemma~\ref{lemma:partition_unbalanced}. Let $v$ be any vertex such that $f_1^+(v) = \min \{f_1^+(x) \mid x\in V(D)\}$. Let $u$ be any vertex such that $f_1^+(u) > f_1^+(v)$. The existence of $u$ is guaranteed, for otherwise the tightness of $(D,F)$ implies that $(D,F)$ is a complete hard pair, a contradiction. Let $X\subseteq (V(D)\setminus \{u,v\})$ be any set of $f_1^+(v)$ vertices (we have $f_1^+(u) \leq |V(D)|-1$, which implies $f_1^+(v) \leq |V(D)|-2$, and the existence of $X$ is guaranteed). Note that $X$, $u$, and $v$ can be computed in linear time.  We colour all the vertices of $X$ with $1$, and note that this is an $F$-dicolouring of $D\ind{X}$, as every vertex $x\in X$ satisfies 
    $f_1^+(x)\geq f_1^+(v)= |X| > |X| - 1 = d_{D\ind{X}}^+(x)$. 
    Let $(D',F'=(f_1',f_2'))$ be the pair reduced from this colouring. Observe that we may have $X=\emptyset$ in which case $(D',F') = (D,F)$. Then $f_1'(v) = (0,0)$ and $f_1'(u) \neq (0,0)$, so $(D',F')$ is not a hard pair. If $(D',F')$ is loose, we are done by Lemma~\ref{lemma:partition_loose}, otherwise we are done by Lemma~\ref{lemma:partition_deficiency} (since $(D',F')$ does not fulfill~\eqref{prop:deficiency}).
\end{proof}

The following shows how to compute, in linear time, an $F$-dicolouring of $D$ if the underlying graph of $D$ is a cycle.
\begin{lemma}
    \label{lemma:partition_cycle}
    Let $(D,F=(f_1,f_2))$ be a valid tight non-hard pair such that $\UG(D)$ is a cycle.
    Then $D$ is $F$-dicolourable and there is a linear-time algorithm providing an $F$-dicolouring of $D$.
\end{lemma}
\begin{proof}
    We assume that~\eqref{prop:deficiency} holds as otherwise we are done by Lemma~\ref{lemma:partition_deficiency}. 
    Let us show that for any vertex $x\in V(D)$, we have $d^-(x),d^+(x)\in \{0,2\}$. Towards a contradiction, and by symmetry, assume that some vertex $x$ verifies $d^-(x)=1$.
    By tightness and by~\eqref{prop:deficiency}, we have that $1= d^-(x)= f_1^-(v) + f_2^-(v) \geq 1+1$, a contradiction.
    We distinguish two cases, depending on whether $D$ contains a digon or not.
    \begin{description}
        \item[Case 1:] \textit{$D$ contains a digon.} 
        
        In this case, there is only one orientation avoiding in- or out-degree one, the one with $D$ fully bidirected. 
        We now claim that, by~\eqref{prop:deficiency-prime}, for every vertex $x\in V(D)$ and every colour $c$, $f_c(x) \geq (1,1)$.
        Hence, by tightness, $f_1,f_2$ are constant functions equal to $(1,1)$. Since $(D,F)$ is not a hard pair, we conclude that $|V(D)|$ is even, and any proper $2$-colouring of $\UG(D)$ is indeed an $F$-dicolouring.
        
        \item[Case 2:] \textit{$D$ does not contain any digon.} 
        
        In this case, the only possible orientation avoiding 
        in- and out-degree one is the antidirected cycle, that is an orientation of a cycle in which every vertex is either a source or a sink. In that case we colour every vertex with colour one, and note that for a sink $x$ (resp. a source) we have that $f_1^+(x)=1>0=d^+(x)$ (resp. $f_1^-(x)=1>0=d^-(x)$). Hence, this colouring is indeed an $F$-dicolouring of $D$.
        \qedhere
    \end{description}
\end{proof}

The following shows how to compute, in linear time, an $F$-dicolouring of $D$ if the underlying graph of $D$ has a cycle going through every vertex but one. We call such a graph a \textit{subwheel}.
\begin{lemma}
    \label{lemma:partition_wheel}
    Let $(D,F=(f_1,f_2))$ be a valid tight non-hard pair, and let $v\in V$ be such that $\UG(D-v)$ is a cycle $C$.
    Then $D$ is $F$-dicolourable and there is a linear-time algorithm providing an $F$-dicolouring of $D$.
\end{lemma}
\begin{proof}
    We assume that~\eqref{prop:deficiency} and~\eqref{prop:digon-sym} hold as otherwise we are done by Lemma~\ref{lemma:partition_deficiency} or Lemma~\ref{lemma:partition_unbalanced}. 
    Let $n=|V(D)|$ and $v_1,\ldots,v_{n-1}$ be an ordering of $V(D)\setminus \{v\}$ along $C$, which can be obtained in linear time. We distinguish several cases, according to the structure of $(D,F)$. It is straightforward to check, in linear time, which of the following cases matches the structure of $(D,F)$.
    
    \begin{description}
        \item[Case 1:] \textit{$|N(v)| < n-1$.}
        
        Let $w$ be any neighbour of $v$ and let $y$ be any vertex that is not adjacent to $v$. Let $\rho_w$ be  $(f_1^-(w) - \mathbb{1}_{A}(v,w),f_1^+(w) - \mathbb{1}_{A}(w,v))$, informally, $\rho_w$ corresponds to the new value of $f_1(w)$ if we were to colour $v$ with $1$. If $\rho_w \neq (1,1)$, we colour $v$ with $1$, otherwise we colour $v$ with $2$. Let $(D',F'=(f_1',f_2'))$ be the pair reduced from colouring $v$, we claim that $(D',F')$ is not a hard pair.
        
        For every $c\in \{1,2\}$, we have $f_c'(y) = f_c(y)$ because $y$ is not adjacent to $v$, and by~\eqref{prop:deficiency} $f'_c(y)=f_c(y) \neq (0,0)$. Therefore $(D',F')$ cannot be a monochromatic hard pair. It is also not a join hard pair because $D'$ is biconnected. 
        
        Assume that $(D',F')$ is a bicycle hard pair. Hence, $f_1'(w) = (1,1)$, so $v$ is coloured $2$, as otherwise we would have $f_1'(w) = \rho_w \neq (1,1)$.
        But if $v$ is coloured $2$, then $(1,1)=f_1'(w) = f_1(w)$ and in that case $\rho_w \neq (1,1)$ (as $w$ is adjacent to $v$), and we should have coloured $v$ with colour 1, a contradiction. 
        
        Assume that $(D',F')$ is a complete hard pair, then $D'$ is a bidirected complete graph. Since $\UG(D')$ is a cycle, $D'$ is necessarily $\bid{K}_3$. Such a complete hard pair is also, either a monochromatic hard pair or a bicycle hard pair, but we already discarded both cases.
        Hence $(D',F')$ is not a hard pair, and the result follows from Lemma~\ref{lemma:partition_cycle}, as $\UG(D')$ is a cycle.
        
        \item[Case 2:] \textit{$D$ is bidirected, $d^+(v) = n-1$, $n$ is even, and $f_1$ is constant on $V(D)\setminus \{v\}$.}
        
        Assume $D\neq \bid{K}_4$ as otherwise the result follows from Lemma~\ref{lemma:partition_complete}.
        Observe that $D-v$ is a bidirected odd cycle, and that it contains at least $5$ vertices for otherwise $D$ is exactly $\bid{K}_4$. Since $D$ is bidirected, by~\eqref{prop:digon-sym}  and~\eqref{prop:deficiency} we have $f_c^+(x) = f_c^-(x) > 0$ for every $x\in V$ and every $c\in \{1,2\}$. 
        
        Since $f_1$ is constant on $V(D)\setminus \{v\}$ by assumption, so is $f_2$, as $(D,F)$ is tight and as $d^-(u)=d^+(u)=3$ for every $u\in V(D)\setminus \{v\}$. Assume without loss of generality that $f_{1}(u)=(1,1)$ and $f_{2}(u)=(2,2)$ for every $u\in V(D)\setminus \{v\}$.
        
        If $f_{1}^+(v) \geq 2$, we colour $v$ and $v_1$ with $1$, and all other vertices with $2$, see Figure~\ref{fig:partition_bid_wheel}$(a)$. The digraph induced by the vertices coloured $2$ is a bidirected path, which is strictly-$f_2$-bidegenerate since on these vertices, $f_2$ is the constant function equal to $(2,2)$. The digraph induced by the vertices coloured $1$ is $\bid{K}_2$ and contains $v$. Since $f_{1}^+(v) \geq 2$ and $f_{1}(v_1) = (1,1)$, it is strictly-$f_{1}$-bidegenerate.

        Else, we have $f_{1}^+(v) \leq 1$, and we colour $\{v\} \cup \{v_{2i} \mid i\in \left[\frac{n}{2}\right] \}$ with $2$ and $\{v_{2i-1} \mid i\in \left[\frac{n-2}{2}\right] \}$ with $1$, see Figure~\ref{fig:partition_bid_wheel}$(b)$.
        Vertices coloured $1$ form an independent set, so the digraph induced by them is strictly-$f_{1}$-bidegenerate (as $f_{1}$ is constant equal to $(1,1)$ on these vertices). Since $f_{1}^+(v) \leq 1$, and because $(D,F)$ is valid, we have $f_{2}^+(v) \geq n-2$. Since $n \geq 6$ it implies $f_2^+(v) \geq \frac{n}{2}+1$. Let $H$ be the digraph induced by the vertices coloured $2$. Then the out-degree of $v$ in $H$ is exactly $\frac{n}{2}$.  Since $f_2^+(v) \geq \frac{n}{2}+1$, we obtain that the digraph $H$ is strictly-$f_2$-bidegenerate if and only if $H-v$ is strictly-$f_2$-bidegenerate.
        Observe that $H-v$ is made of a copy of $\bid{K}_2$ and isolated vertices, so it is strictly-$f_2$-bidegenerate since $f_2$ is constant equal to $(2,2)$ on $V(H)\setminus \{v\}$.
        
        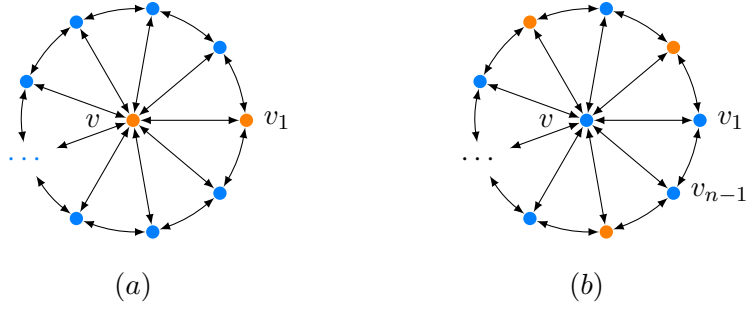
\begin{figure}[hbtp]
            \begin{minipage}{\linewidth}
                \begin{center}	
                    \begin{tikzpicture}[scale=1, every node/.style={transform shape}]
                        \tikzset{vertex/.style = {circle,fill=black,minimum size=5pt, inner sep=0pt}}
                        \tikzset{edge/.style = {->,> = latex'}}
                        \tikzset{digon/.style = {<->,> = latex}}
                        
                        \node[vertex, orange, label={[label distance=0.2cm]left:$v$}] (v) at (0,0) {};
                        \node[] (a) at (0,-2.2) {$(a)$};
                        \node[vertex, orange, label=right:$v_1$] (v0) at (0:1.5) {};
                        \node[g-blue] (v5) at (5*360/9:1.5) {$\cdots$};
                        \foreach \i in {1,2,3,4,6,7,8}{
                            \node[vertex,g-blue] (v\i) at (\i*360/9:1.5) {};
                        }
                        \foreach \i in {0,...,8}{
                            \pgfmathtruncatemacro{\j}{Mod(\i+1,9)}
                            \draw[digon, bend right=10] (v\i) to (v\j) {};
                            \draw[digon] (v\i) to (v) {};
                        }
                        
                        \begin{scope}[xshift=6cm]
                            \node[vertex, g-blue, label={[label distance=0.2cm]left:$v$}] (v) at (0,0) {};
                            \node[] (b) at (0,-2.2) {$(b)$};
                            \node[vertex, g-blue, label=right:$v_1$] (v0) at (0:1.5) {};
                            \node[] (v5) at (5*360/9:1.5) {$\cdots$};
                            \node[vertex,g-blue,label=right:$v_{n-1}$] (v8) at (8*360/9:1.5) {};
                            \foreach \i in {1,3,7}{
                                \node[vertex,orange] (v\i) at (\i*360/9:1.5) {};
                            }
                            \foreach \i in {2,4,6}{
                                \node[vertex,g-blue] (v\i) at (\i*360/9:1.5) {};
                            }
                            \foreach \i in {0,...,8}{
                                \pgfmathtruncatemacro{\j}{Mod(\i+1,9)}
                                \draw[digon, bend right=10] (v\i) to (v\j) {};
                                \draw[digon] (v\i) to (v) {};
                            }
                        \end{scope}
                    \end{tikzpicture}
                    \caption{The $F$-dicolourings for Case~2. The partition on the left corresponds to the case $f^+_{1}(v)\geq 2$ and the one on the right corresponds to the case $f^+_{1}(v)\leq 1$. Vertices coloured $1$ are represented in orange and vertices coloured $2$ are represented in blue.}
                    \label{fig:partition_bid_wheel}
                \end{center}    
            \end{minipage}
        \end{figure}
        
        \item[Case 3:] \textit{$N^+(v) = N^-(v) = V(D)\setminus \{v\}$ and $D-v$ is a directed cycle}. 
        
        Note that in this case, for every vertex $u\in V(D)\setminus \{v\}$ we have $d^-(u)=d^+(u)=2$, and by~\eqref{prop:deficiency} and by tightness of $(D,F)$, we have that $f_1(u)=f_2(u)=(1,1)$.
        Since $d^+(v) \geq 3$, there is a colour $c\in \{1,2\}$ such that $f_c^+(v)\geq 2$. Assume $f_2^+(v)\geq 2$. 
        We are now ready to give an $F$-dicolouring explicitly. 
        We colour $v$ and $v_1$ with $2$, and all other vertices with $1$, see Figure~\ref{fig:partition_simple_wheel} for an illustration.
        
        The digraph induced by the vertices coloured $1$ is a directed path, which is strictly-$f_{1}$-bidegenerate since on these vertices, $f_{1}$ is the constant function equal to $(1,1)$. The digraph induced by the vertices coloured $2$ is $\bid{K}_2$ and contains $v$. Since $f_{2}^+(v) \geq 2$ and $f_{2}(v_1) = (1,1)$, it is strictly-$f_{2}$-bidegenerate.
        \begin{figure}[hbtp]
            \begin{minipage}{\linewidth}
                \begin{center}	
                    \begin{tikzpicture}[thin,scale=1, every node/.style={transform shape}]
                        \tikzset{vertex/.style = {circle,fill=black,minimum size=5pt, inner sep=0pt}}
                        \tikzset{edge/.style = {->,> = latex}}
                        \tikzset{digon/.style = {<->,> = latex}}
                        
                        \node[vertex, g-blue, label={[label distance=0.2cm]left:$v$}] (v) at (0,0) {};
                        \node[vertex, g-blue, label=right:$v_1$] (v0) at (0:1.5) {};
                        \node[orange] (v5) at (5*360/7:1.5) {$\cdots$};
                        \foreach \i in {1,2,3,4,6}{
                            \node[vertex,orange] (v\i) at (\i*360/7:1.5) {};
                        }
                        \foreach \i in {0,...,6}{
                            \pgfmathtruncatemacro{\j}{Mod(\i+1,7)}
                            \draw[edge, bend right=15] (v\i) to (v\j) {};
                            \draw[digon] (v\i) to (v) {};
                        }
                    \end{tikzpicture}
                    \caption{The $F$-dicolouring for Case~3. Vertices coloured $1$ are represented in orange and vertices coloured $2$ are represented in blue.}
                    \label{fig:partition_simple_wheel}
                \end{center}    
            \end{minipage}
        \end{figure}
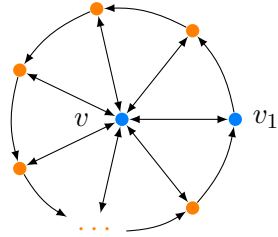
        
        \item[Case 4:] \textit{None of the previous cases apply.}
        
        We first prove the existence of a vertex $u\neq v$ and a colour $c\in \{1,2\}$ such that $f_c(u) > (\mathbb{1}_A(v,u), \mathbb{1}_A(u,v))$. By~\eqref{prop:deficiency-prime}, for any vertex $u\neq v$ we have $f_c(u) \geq (\mathbb{1}_A(v,u), \mathbb{1}_A(u,v))$, we thus look for a vertex $u\neq v$ such that        
        $f_c(u) \neq (\mathbb{1}_A(v,u), \mathbb{1}_A(u,v))$. Assume for a contradiction that such a pair does not exist, so for every vertex $x\in V(D)\setminus \{v\}$ and every colour $c\in \{1,2\}$, $f_c(x)$ is exactly $(\mathbb{1}_{A}(v,x),\mathbb{1}_{A}(x,v))$. Assume first that there exists a simple arc $xv$, then $f_1(x) = f_2(x) = (0,1)$. By~\eqref{prop:deficiency}, we deduce that $x$ is a source. Since $\UG(D-v)$ is a cycle, we then have $d^+(x) = 3 > f_1^+(x) + f_2^+(x)$, a contradiction to $(D,F)$ being a valid pair. The existence of a simple arc $vx$ is ruled out symmetrically, so we now assume $N^+(v) = N^-(v)$. Then every vertex $x\in V(D)\setminus \{v\}$ satisfies $f_1(x) = f_2(x) = (\mathbb{1}_{A}(v,x),\mathbb{1}_{A}(x,v)) = (1,1)$.
        Since $(D,F)$ is a tight valid pair, this implies that $d^+(x) = d^-(x) = 2$. Hence $D-v$ is a directed cycle, so Case~3 is matched, a contradiction. This proves the existence of $u$ and $c$ such that $f_c(u) \neq (\mathbb{1}_A(v,u), \mathbb{1}_A(u,v))$.
        From now on, we assume $c=1$ without loss of generality, and $f_1(u) > (\mathbb{1}_A(v,u), \mathbb{1}_A(u,v))$.

        Consider the following property, which can be checked straightforwardly in linear time:
        \begin{equation}
            \exists x \in V(D)\setminus \{v\},~~\text{such that}~~f_1(x) \neq (1+\mathbb{1}_{A}(v,x), 1+\mathbb{1}_{A}(x,v)) \text{~~~or~~~} f_{2}(x) \neq (1, 1).\tag{3}\label{eq:wheel:property_x}
        \end{equation}
        If $n$ is odd or if~\eqref{eq:wheel:property_x} holds, we colour $v$ with $1$. Otherwise, we colour $v$ with $2$. In both cases, we claim that the reduced pair $(D',F'=(f_1',f_2'))$ is not hard. 
        
        Assume first that $n$ is odd or~\eqref{eq:wheel:property_x} is satisfied, so we have coloured $v$ with $1$. 
        Hence, $f_{1}'(u) = f_{1}(u) - (\mathbb{1}_A(v,u), \mathbb{1}_A(u,v)) \neq (0,0)$ and $f_{2}'(u) = f_{2}(u) \neq (0,0)$ (by~\eqref{prop:deficiency}), so $(D',F')$ is not a monochromatic hard pair. If $n$ is odd, then $|V(D')|$ is even so $D'$ is not a bidirected odd cycle. If a vertex $x$ satisfies~\eqref{eq:wheel:property_x}, we have $f_{1}'(x) = f_{1}(x) - (\mathbb{1}_A(v,x), \mathbb{1}_A(x,v)) \neq (1,1)$ or $f_{2}'(x)=f_2(x) \neq (1,1)$. In both cases $(D',F')$ is not a bicycle hard pair, as desired. Since $UG(D')$ is a cycle, if $(D',F')$ was a complete hard pair it would be either a monochromatic or a bicycle hard pair, hence $(D',F')$ is not a hard pair.
        
        Assume finally that $n$ is even and~\eqref{eq:wheel:property_x} is not satisfied. Hence $v$ is coloured with $2$ and, for every vertex $x\in V(D)\setminus \{v\}$, $f_1(x) = (1+\mathbb{1}_{A}(v,x), 1+\mathbb{1}_{A}(x,v))$ and $f_{2}(x) = (1,1)$. Since every vertex $x\neq v$ is adjacent to $v$, otherwise Case~$1$ would match, we have $f'_2(x) \neq f_2(x) = (1,1)$ so $(D',F')$ is not a bicycle hard pair. Assume for a contradiction that it is a monochromatic hard pair. Since $f'_1(x)=f_1(x)\neq (0,0)$ for every vertex $x\neq v$, we have $f_{2}'(x) = (0,0)$. Since $f_{2}(x)=(1,1)$, we deduce that $N^+(v) = N^-(v) = V(D) \setminus \{v\}$. Hence every vertex $x$ satisfies $f_1(x) = (2,2)$ and $f_{2}=(1,1)$. Since $(D,F)$ is tight, we conclude that $D$ is bidirected, $n$ is even and $f_1$ is constant on $V(D) \setminus \{v\}$. Thus Case~2 matches, a contradiction.
        
        Since $(D',F')$ is not a hard pair, and because $\UG(D')$ is a cycle, the result follows from Lemma~\ref{lemma:partition_cycle}.
        \qedhere
    \end{description}
\end{proof}

\subsection{Solving blocks with two colours -  general case}\label{subsec:algo-block-2col-general}

The goal of this subsection is to provide an algorithm which, given a pair $(D,F=(f_1,f_2))$ that is not hard and such that $D$ is biconnected, computes an $F$-dicolouring of $D$ in linear time.
The idea is to find a decomposition of the underlying graph $\UG(D)$, similar to an ear-decomposition, and successively colour the ``ears'' in such a way that the successive reduced pairs remain non-hard.
Let us begin with the definition of the decomposition.
\begin{definition}\label{def:CSP-decomp}
    A \textit{CSP-decomposition} (CSP stands for Cycle, Stars, and Paths) of a biconnected graph $G$ is a sequence $(H_0,\ldots, H_r)$ of subgraphs of $G$ partitioning the edges of $G$, such that $H_0$ is a cycle, and such that for any $i\in [r]$ the subgraph $H_i$ is either:
    \begin{itemize}
        \item a star with a central vertex $v_i$ of degree at least two in $H_i$, and such that $V(H_i) \cap \left(\bigcup_{0\leq j< i} V(H_j) \right)$ is the set of leaves of $H_i$, or
        \item a path $(v_0,\ldots,v_\ell)$ of length $\ell\geq 3$, and such that $V(H_i) \cap \left(\bigcup_{0\leq j< i} V(H_j) \right) = \{v_0,v_\ell\}$. 
    \end{itemize}
\end{definition}
\begin{lemma}\label{lemma:CSP-decomp}
    Every biconnected graph admits a CSP-decomposition. Furthermore, computing such decomposition can be done in linear time.
\end{lemma}
\begin{proof}
    It is well-known that every biconnected graph $G$ admits an ear-decomposition, and that it can be computed in linear time (see~\cite{schmidtIPL113} and the references therein).
    An \textit{ear-decomposition} is similar to a CSP-decomposition except that paths of length one or two are allowed, and that there are no stars. To obtain a CSP-decomposition of $G$, we first compute an ear-decomposition $(H_0,\ldots, H_r)$, which we modify in order to get a CSP-decomposition.
    We do this in two steps.
    
    Before describing these steps, we have to explain how a decomposition $(H_0,\ldots, H_r)$ is encoded. The sequence is encoded as a doubly-linked chain whose cells contain 1) a copy of the subgraph $H_i$ 2) a mapping from the vertices of this copy to their corresponding ones in $G$, and 3) an integer $n_i$ equal to $|\bigcup_{0\leq j< i} V(H_j)|$. This last integer will allow us, given two cells, corresponding to $H_i$ and $H_j$, to decide whether $i\leq j$ or not. Indeed, as the decomposition will evolve along the execution, maintaining indices from $\{0,\ldots,r\}$ seems hard in linear time.
    
    The first step consists in modifying the ear-decomposition in order to get an ear-decomposition avoiding triples $i\leq j <k$ such that $H_k=(u_0,u_1)$ is of length one, where $u_0,u_1$ are inner vertices of $H_i, H_j$ respectively, and where $H_j$ is an ear of length at least three.
    Towards this, we first go along every ear, in increasing order, to compute $\texttt{birth}(v)$, the copy of vertex $v$ appearing first in the current ear-decomposition.
    Note that $\texttt{birth}(v)$ is always an inner vertex of its ear. Then, we go along every ear in decreasing order, and when we have a length one ear $H_k = (u_0,u_1)$, we consider the ears of $\texttt{birth}(u_0)$ and $\texttt{birth}(u_1)$, say $H_i$ and $H_j$ respectively.
    
    If $0\leq i<j$, and $H_j$ is a path of length at least three, we combine $H_j$ and $H_k$ into two ears, $H'$ and $H''$, each of length at least two (see Figure~\ref{fig:ears} (left)). In the ear-decomposition, we insert $H'$ and $H''$ at the position of $H_j$, and delete $H_j$ and $H_k$.
    
    If $0<i=j$, then $H_j=H_i$ must be a path of length at least four. In that case, we combine $H_j$ and $H_k$ into two ears, $H'$ and $H''$, each of length at least two (see Figure~\ref{fig:ears} (middle)). In the ear-decomposition, we insert $H'$ and $H''$ at the position of $H_j$, and delete $H_j$ and $H_k$.
    
    If $0=i=j$, then $H_j=H_0$ must be a cycle of length at least four. In that case, we combine $H_j$ and $H_k$ into two ears, $H'$ and $H''$, each of length at least two (see Figure~\ref{fig:ears} (right)). In the ear-decomposition, we insert $H'$ and $H''$ at the position of $H_j$, and delete $H_j$ and $H_k$.
    
    \begin{figure}[hbtp]
        \begin{center}
            \includegraphics[width=.95\textwidth]{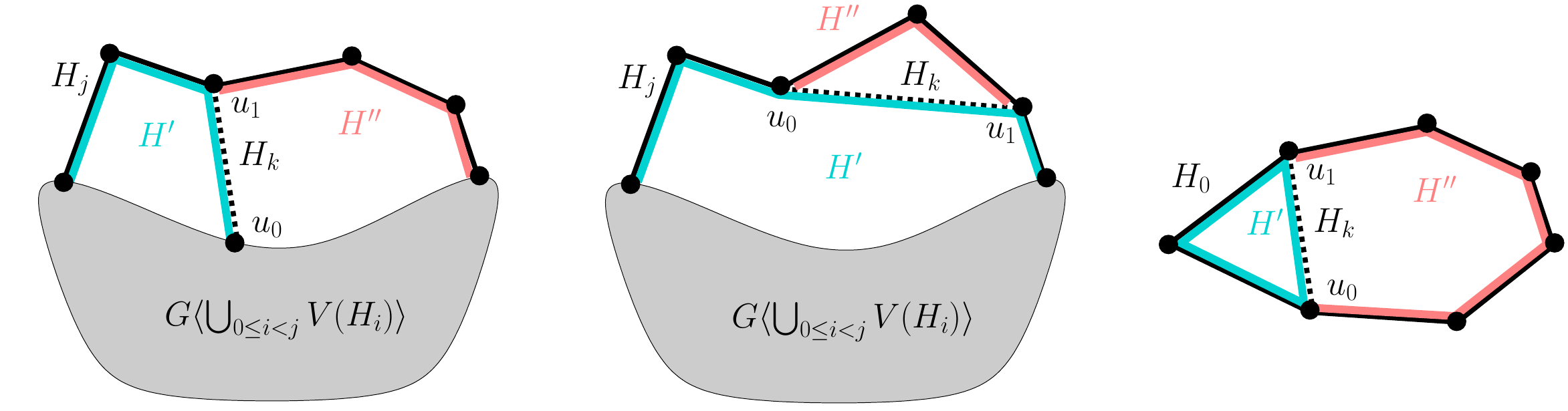} 
            \caption{Ears $H_j$ and $H_k$ recombined into $H'$ (light blue) and $H''$ (pink), each of length at least two. On the left, the case where $H_k$ has one endpoint in $H_j$. In the middle, the case where $H_k$ has its two endpoints in $H_j$. On the right, the case where $H_k$ has its two endpoints in the cycle $H_j=H_0$.}
            \label{fig:ears}
        \end{center}
    \end{figure}
    
    Along this process, the ears that are modified are always shortened while still having length at least two. Hence, we do not create new couples $H_j,H_k$, with the second of length one, that could be recombined into ears of length at least two. Hence, we are done with the first step.
    
    \medskip
    
    Before proceeding to the second step, recall that we now have an ear-decomposition without triples $i\leq j <k$ such that $H_k=(u_0,u_1)$ is of length one, where $u_0,u_1$ are inner vertices of $H_i, H_j$ respectively, and where $H_j$ is an ear of length at least three. In particular, this implies that for any ear of length one $H_k=(u_0,u_1)$, the ears containing $u_0$ and $u_1$ as inner vertices are distinct.
    We denote them $H_i$ and $H_j$, respectively, and assume without loss of generality that $i<j$. Note that as $H_j$ and $H_k$ cannot be recombined, we have that $H_j$ has length exactly two and that $u_1$ is its only inner vertex. The second step then simply consists, for each such ear $H_k$, in including $H_k$ into $H_j$, redefining the latter to form a star centred in $u_1$.
    
    After this process, no ear can be a path of length one. Assimilating the paths of length two as stars, we thus have a CSP-decomposition. 
\end{proof}

We move to showing that $(D,F)$ may be inductively coloured following this decomposition (on the underlying graph of $D$).
Again, we only consider instances that are not dealt with in the previous subsection, which in particular satisfy~\eqref{prop:deficiency} and~\eqref{prop:digon-sym}.
Then, in the induction, Lemma~\ref{lemma:treat_single_ear} corresponds to colouring a star, and Lemma~\ref{lemma:treat_large_ear} corresponds to colouring a path, when those attach to the rest of the graph as in the decomposition. Eventually, bidirected complete graphs, cycles, and subwheels, which were solved in the previous subsection, form the base of our decomposition. The induction is formalised in Lemma~\ref{lemma:solving_2colours}, which achieves to solve the case of biconnected graphs with two colours. 

\begin{lemma}
    \label{lemma:treat_single_ear}
    Let $(D=(V,A),F=(f_1,f_2))$ be a valid tight non-hard pair, and let $v\in V(D)$ be such that both $D$ and $D-v$ are biconnected, and $\UG(D)$ is neither a cycle nor a subwheel centred in $v$. 
    Then there exists a colour $c\in \{1,2\}$ such that colouring $v$ with $c$ is safe. Moreover, there is an algorithm that either finds $c$ in time $O(d(v))$ or computes an $F$-dicolouring of $D$ in linear time.
\end{lemma}
\begin{proof}
    We assume that~\eqref{prop:deficiency} and~\eqref{prop:digon-sym} hold as otherwise we are done by Lemma~\ref{lemma:partition_deficiency} and Lemma~\ref{lemma:partition_unbalanced}. 
    Observe that $d(v) > 0$ and~\eqref{prop:deficiency} implies that $f_1(v) \neq (0,0)$ and $f_2(v) \neq (0,0)$, so both colours $1$ and $2$ are available for $v$. 
    Note also that since $\UG(D-v)$ is not a cycle, in particular colouring $v$ with any colour $c\in \{1,2\}$ may never reduce it to a bicycle hard pair. We distinguish three cases.
    
    \begin{description}
        \item[Case 1:] \textit{$|N(v)| \leq |V(D)|-2$.} 
        
        Let $u$ be any neighbour of $v$ and $w$ be any vertex that is not adjacent to $v$. Since $|N(w)| > 0$, we have $f_1(w) \neq (0,0)$ and $f_2(w) \neq (0,0)$ by~\eqref{prop:deficiency}. Hence colouring $v$ with any colour $c\in \{1,2\}$ does not reduce to a monochromatic hard pair. If $f_1(u)$ is distinct from $(f_1^-(w) + \mathbb{1}_{A}(vu),f_1^+(w) + \mathbb{1}_{A}(uv))$, we colour $v$ with $1$, otherwise we colour $v$ with $2$. By choice of $u$ and $w$, the reduced pair $(D',F'=(f_1',f_2'))$ satisfies $f_1'(u) \neq f_1'(w)$, so in particular it is not a complete hard pair. 
        
        \item[Case 2:] \textit{$|N(v)| = |V(D)|-1$ and $D-v$ is distinct from $\bid{K}_{|V(D)|-1}$.} 
        
        We only have to guarantee the existence of $c\in \{1,2\}$ such that colouring $v$ with $c$ does not reduce to a monochromatic hard pair. Also note that $|N(v)|=|V(D)|-1$ allows us $O(|V(D)|)$ time for computing $c$.
        
        If there exists a vertex $x\neq v$ and a colour $c\in \{1,2\}$ such that $f_c(x) > (\mathbb{1}_A(vx),\mathbb{1}_A(xv))$, then colouring $v$ with $c$ is safe as it does not reduce to a monochromatic hard pair as $f'_c(x)= f_c(x) - (\mathbb{1}_A(vx),\mathbb{1}_A(xv)) \neq (0,0)$ and as $f'_{c'}(x)= f_{c'}(x)\neq (0,0)$ by~\eqref{prop:deficiency} for the colour $c'\neq c$. Note that, if it exists, such a vertex $x$ can be found in time $O(|V(D)|)$.
        
        We now prove the existence of such a vertex $x$. Assume for a contradiction that $f_c(x) \ngtr (\mathbb{1}_A(vx),\mathbb{1}_A(xv))$ for every $x\neq v$ and every $c\in \{1,2\}$.
        Observe that, for every $x\neq v$, $f_c(x)\geq (\mathbb{1}_A(vx),\mathbb{1}_A(xv))$ by~\eqref{prop:deficiency-prime}, so we have $f_1(x) = f_2(x) = (\mathbb{1}_A(vx),\mathbb{1}_A(xv))$. Let $x\neq v$ be any vertex. If $xv$ is a simple arc, then $f_1(x) = f_2(x) = (0,1)$ so $d^+(x) = 2$ and $d^-(x)=0$. In particular, $x$ has degree $1$ in $\UG(D-v)$, a contradiction since $\UG(D-v)$ is biconnected and distinct from $K_2$ (as $\UG(D)$ is not a cycle). The case of $vx$ being a simple arc is symmetric, so we now assume $N^+(v) = N^-(v) = V \setminus \{v\}$, which implies $f_1(x) = f_2(x) = (1,1)$ for every vertex $x\neq v$. 
        Hence in $D-v$ every vertex has in- and out-degree one. Since $\UG(D-v)$ is biconnected and distinct from $K_2$, $D-v$ is necessarily a directed cycle, a contradiction.
        
        \item[Case 3:] \textit{$|N(v)| = |V(D)|-1$ and $D-v$ is $\bid{K}_{|V(D)|-1}$.}
        
        If $N^+(v) = N^-(v)$ then $D$ is $\bid{K}_{|V(D)|}$, and since $(D,F)$ is not hard, Lemma~\ref{lemma:partition_complete} yields that it is $F$-dicolourable and an $F$-dicolouring is computed in linear time.
        
        Henceforth we assume that there exists a simple arc between $v$ and some $u\in V(D) \setminus \{v\}$.
        The vertex $u$ being incident to some digon,~\eqref{prop:deficiency} implies that $f_c(u)\geq (1,1)$ for every colour $c$. This guarantees that colouring $v$ with any colour $c\in \{1,2\}$ does not reduce to a monochromatic hard pair, since in the reduced pair the constraints for $u$ are unchanged on one coordinate. We thus have to guarantee that for some $c\in \{1,2\}$, colouring $v$ with $c$ does not reduce to a complete hard pair.
        
        We assume that the simple arc between $v$ and $u$ goes from $u$ to $v$, the other case being symmetric. If $v$ has at least one out-neighbour $w$, then either $f_1^-(u)=f_1^-(w)$ and colouring $v$ with $1$ does not reduce to a complete hard pair, or $f_1^-(u)\neq f_1^-(w)$ and colouring $v$ with $2$ does not reduce to a complete hard pair. 
        
        Finally, if $v$ is a sink and colouring $v$ with $1$ reduces to a complete hard pair, then $f_1^+(x) = f_1^-(x) + 1$ for every vertex $x\neq v$. Therefore, colouring $v$ with $2$ does not reduce to a complete hard pair, as in the reduced pair we have $f_1'(x) = f_1(x)$ and $f_1^+(x) \neq f_1^-(x)$.
        \qedhere
    \end{description}
\end{proof}

\begin{lemma}
    \label{lemma:treat_large_ear}
    Let $(D=(V,A),F=(f_1,f_2))$ be a valid pair satisfying~\eqref{prop:deficiency} and let $v_0,\ldots,v_{\ell}$ be a path of $G = \UG(D)$ of length $\ell\geq 3$ such that $d_G(v_i)=2$ for $i\in[\ell-1]$, and such that both $G = \UG(D)$ and $G' = \UG(D - \{v_1,\ldots,v_{\ell-1}\})$ are biconnected and contain a cycle.
    There is an algorithm computing an $F$-dicolouring of $D[\{v_1,\ldots,v_{\ell-1}\}]$ in time $O(\ell)$ in such a way that the reduced pair $(D',F')$ is not hard.
\end{lemma}
\begin{proof}
    Since~\eqref{prop:deficiency} is satisfied, in particular, for any vertex $x$ and colour $c$, we have $f_c(x)\neq (0,0)$.  
    Note that $V \setminus \{v_0,\ldots,v_\ell\} \neq \emptyset$ for otherwise $G'$ does not contain a cycle, a contradiction.
    Let $u$ be any vertex in $V \setminus \{v_0,\ldots,v_\ell\}$. If $f_1(v_0) = f_1(u)$, we colour $v_1$ with $1$, otherwise we colour $v_1$ with $2$. We then (greedily) colour $v_2,\ldots,v_{\ell-1}$ using the algorithm of Lemma~\ref{lemma:validity-algo-greedy}.
    For every $j\in \{2,\ldots,\ell-1\}$, the vertex $v_j$ has at least one neighbour in $\{v_{j+1},\ldots,v_{\ell}\}$, $v_{j+1}$,  so condition~\eqref{property_degenerate_F} is fulfilled, and we are ensured that the obtained colouring is an $F$-dicolouring $D\ind{\{v_1,\ldots,v_{\ell-1}\}}$.
    
    Let $(D',F'=(f_1',f_2'))$ be the pair reduced from this colouring. 
    As $D'$ is biconnected, $(D',F')$ is not a join hard pair.
    By construction, we have $f_1'(v_0) \neq f_1'(u)$ as $v_0$ is adjacent to $v_1$ but $u$ is not adjacent to any $v_i$ with $i \in [1,\ell-1]$, which we just coloured. Hence $(D',F')$ is neither a bicycle hard pair nor a complete hard pair, as $f_1'$ is not constant on $V(D')$. Finally,~\eqref{prop:deficiency} ensures that $f_1'(u)=f_1(u) \neq (0,0)$ and $f_2'(u) = f_2(u) \neq (0,0)$. Hence $(D',F')$ is not a monochromatic hard pair.
\end{proof}

With these lemmas in hand we are ready to prove the main result of this subsection, that non-hard biconnected pairs are $F$-dicolourable.
\begin{lemma}
    \label{lemma:solving_2colours}
    Let $(D=(V,A),F=(f_1,f_2))$ be a valid non-hard pair. There exists an algorithm providing an $F$-dicolouring of $D$ in linear time.
\end{lemma}
\begin{proof}
    We denote by $\mathcal{P}$ the property of $(D,F)$ being tight, fulfilling~\eqref{prop:deficiency} and~\eqref{prop:digon-sym}, not being a bidirected complete graph.
    
    We first check in linear time that each property of $\mathcal{P}$ holds.
    If one does not, we are done by one of Lemmas~\ref{lemma:partition_loose},~\ref{lemma:partition_deficiency},~\ref{lemma:partition_unbalanced} and~\ref{lemma:partition_complete}.
    
    We then compute a CSP-decomposition $(H_0,\ldots,H_r)$ of $\UG(D)$ in linear time, which is possible by Lemma~\ref{lemma:CSP-decomp}. If $r=0$ then $\UG(D)$ is a cycle and the result follows from Lemma~\ref{lemma:partition_cycle}, assume now that $r\geq 1$.
    
    For each $i$ going from $r$ to $2$, we proceed as follows (if $r=1$ we skip this part). If $H_i$ is a path of length $\ell \geq 3$, we colour $v_1,\ldots,v_{\ell-1}$ in time $O(\ell)$ in such a way that the reduced pair $(D',F')$ is not hard, which is possible by Lemma~\ref{lemma:treat_large_ear}. We then check in constant time that $(D',F')$ satisfies $\mathcal{P}$ (we only have to check that it still holds for $v_0$ and $v_\ell$, and this is done in constant time). If it does not, we are done by one of Lemmas~\ref{lemma:partition_loose},~\ref{lemma:partition_deficiency},~\ref{lemma:partition_unbalanced} and~\ref{lemma:partition_complete}.
    If $H_i$ is a star with central vertex $v$, then we compute a safe colouring of $v$ in time $O(d(v))$, which is possible by Lemma~\ref{lemma:treat_single_ear}. Note that, at this step, we may directly find an $F$-dicolouring in linear time, in which case we stop here. Assume we do not stop, then we obtain a reduced pair $(D',F')$. We then check in time $O(d(v))$ that $(D',F')$ satisfies $\mathcal{P}$ (we only have to check that it still holds for the neighbours of $v$ in $H_i$). Again, if it does not, we are done by one of Lemmas~\ref{lemma:partition_loose},~\ref{lemma:partition_deficiency},~\ref{lemma:partition_unbalanced} and~\ref{lemma:partition_complete}.
    
    Assume we have not already found an $F$-dicolouring by the end of this process, and consider $H_1$. If $H_1$ is a star, we are done by Lemma~\ref{lemma:partition_wheel}. Otherwise, $H_1$ is a path, and we colour it as in the process above, then we conclude by colouring $H_0$ through Lemma~\ref{lemma:partition_cycle}.
    
    Since $(H_0,\ldots,H_r)$ partitions the edges of $\UG(D)$, the total running time of the described algorithm is linear in the size of $D$.
\end{proof}

\subsection{Proof of Theorem~\ref{thm:bivariable}}
\label{subsec:proof_thm_bivariable}

We are now ready to prove Theorem~\ref{thm:bivariable}, that we first recall here for convenience.

\bivariable*
\begin{proof}
    Let $(D,F)$ be a valid pair. We first check in linear time whether $(D,F)$ is tight. If it is not, we are done by Lemma~\ref{lemma:partition_loose}. Henceforth assume that $(D,F)$ is tight. 
    We then check whether $(D,F)$ is a hard pair, which is possible in linear time by Lemma~\ref{lemma:check_hard-pair}. If $(D,F)$ is a hard pair, then $D$ is not $F$-dicolourable by Lemma~\ref{lemma:hard_pairs_not_dicolourable} and we can stop.
    Then, we may assume that $(D,F)$ is not a hard pair.
    In linear time, we find a block $B$ of $D$ and compute a safe colouring $\alpha$ of $V(D) \setminus V(B)$, which is possible by Lemma~\ref{lemma:reduction_to_block}.
    Let $(B,F_B)$ be the reduced pair, which is not hard, 
    then every $F_B$-dicolouring of $B$, together with $\alpha$, extends to an $F$-dicolouring of $D$.
    Then, in linear time, we either find an $F_B$-dicolouring of $B$, or we compute a valid pair $(B,\Tilde{F}_B=(\Tilde{f}_1,\Tilde{f}_2))$ such that $(B,\Tilde{F})$ is not hard, and an $F_B$-dicolouring of $B$ can be computed in linear time from any $\Tilde{F}_B$-dicolouring of $B$. This is possible by Lemma~\ref{lemma:reduce_to_2colours}.
    
    If we have not yet found an $F_B$-dicolouring of $B$, we compute an $\Tilde{F}_B$-dicolouring of $B$ in linear time, which is possible by Lemma~\ref{lemma:solving_2colours}. From this we compute an $F_B$-dicolouring of $B$ in linear time. As mentioned before, this $F_B$-dicolouring together with $\alpha$ gives an $F$-dicolouring of $D$, obtained in linear time.
\end{proof}

\end{document}